\documentclass[12pt]{article}
\usepackage{amsmath,amsxtra,amssymb,color,latexsym,epsfig,amscd,amsthm,subfigure,fancybox,epsfig}
\usepackage[mathscr]{eucal}
\usepackage{bbm}
\usepackage{graphicx}
\usepackage{enumerate}
\usepackage{enumitem}
\usepackage{epsfig}
\usepackage{epstopdf}
\usepackage{cases}
\def\para#1{\vskip .4\baselineskip\noindent{\bf #1}}
\def\tr{{\rm tr}}

\newtheorem{thm}{Theorem}[section]
\newtheorem {asp}{Assumption}[section]
\newtheorem{lm}{Lemma}[section]
\newtheorem{prop}{Proposition}[section]

\theoremstyle{definition}
\newtheorem{defn}[thm]{Definition}

\theoremstyle{remark}
\newtheorem{rem}{Remark}[section]
\newtheorem{exam}{Example}[section]
\numberwithin{equation}{section}

\newcommand{\eps}{\varepsilon}

\newcommand{\h}{\mathcal{H}}

\newcommand{\M}{\mathcal{M}}
\newcommand{\F}{\mathcal{F}}

\newcommand{\E}{\mathbb{E}}

\newcommand{\Lom}{\mathcal{L}}

\newcommand{\N}{{\mathbb{Z}}_+}

\newcommand{\PU}{\mathcal{P}(\R^{2,\circ}_+\times[0,M])}
\newcommand{\UU}{\R^{2,\circ}_+\times[0,M]}
\newcommand{\PUB}{\mathcal{P}((0,\infty)\times[0,\infty)\times[0,M])}
\newcommand{\PP}{\mathbb{P}}

\newcommand{\R}{\mathbb{R}}

\numberwithin{equation}{section}
\newcommand{\1}{\boldsymbol{1}}

\newcommand{\blue}[1]{{\leavevmode\color{black}#1}}

\newcommand{\wdt}{\widetilde}

\newcommand{\op}{{\cal L}}

\newcommand{\bed}{\begin{displaymath}}
\newcommand{\eed}{\end{displaymath}}
\newcommand{\bea}{\bed\begin{array}{rl}}
\newcommand{\eea}{\end{array}\eed}
\newcommand{\ad}{&\!\!\!\disp}
\newcommand{\aad}{&\disp}
\newcommand{\barray}{\begin{array}{ll}}
\newcommand{\earray}{\end{array}}
\newcommand{\diag}{{\rm diag}}

\def\disp{\displaystyle}
\newcommand{\rr}{{\Bbb R}}

\def\bar{\overline}
\def\hat{\widehat}
\def\a.s{\text{\;a.s.\;}}

\begin{document}
\title{General Nonlinear Stochastic Systems Motivated by
Chemostat Models:
Complete Characterization of
Long-Time Behavior,
Optimal Controls, and
 Applications to
Wastewater Treatment}
\author{Dang H. Nguyen,\thanks{Department of Mathematics
     University of Alabama
     Tuscaloosa, AL 35401,
dangnh.maths@gmail.com. The research
 of this author was supported in part by the National Science Foundation under grant DMS-1853467.}
\and Nhu N. Nguyen,\thanks{Department of Mathematics, Wayne State University, Detroit, MI
48202,
nhu.math.2611@gmail.com. The research
 of this author was supported in part by the National Science Foundation under grant DMS-1710827.}
\and
George Yin\thanks{Department of Mathematics, Wayne State University, Detroit, MI
48202,
gyin@wayne.edu. The research
 of this author was supported in part by the National Science Foundation under grant DMS-1710827.}}
\date{}
\maketitle

\begin{abstract}
This paper focuses on a general
class of systems of nonlinear stochastic differential equations, inspired by stochastic chemostat models. In the first part,
the system is formulated  as a hybrid switching diffusion.
A complete characterization of the asymptotic behavior of the system under consideration is provided. It is shown that the long-term properties of the system can be classified by using a real-valued parameter $\lambda$. If $\lambda\leq 0$, the bacteria
will die out, which means that the process does not operate.
If $\lambda>0$, the system has an invariant probability measure
and the transition probability of the solution process converges to that of the invariant measure. The rate of convergence is also obtained. One of the distinct features of this paper is that the critical case $\lambda=0$ is also considered. Moreover,
numerical examples are given to illustrate our results.
In the second part of the paper, controlled diffusions
 with a long-run average objective function are treated.
The associated Hamilton-Jacobi-Bellman (HJB) equation is derived and the existence of an optimal Markov control is established.
The techniques and methods of analysis in this paper can be applied to many other stochastic Kolmogorov systems.

\bigskip

\end{abstract}

\newpage
\section{Introduction}\label{sec:int}
We consider a class of systems of nonlinear stochastic differential equations. The motivation stems from the pressing need of the treatment of chemostat models that are a laboratory apparatus used for the continuous culture of microorganism, which is
a technique introduced by  Novick and Szilard in \cite{NS50}. This technique plays an important role in microbiology, biotechnology, and population biology,  and is perhaps the best laboratory idealization of nature for population studies \cite{W71}.
Chemostats are also used as microcosms in ecology
\cite{BHMJA05, PK92} and evolutionary biology \cite{DD04,JE07},
as well as in wastewater treatment based on chemostat models
\cite{BK, GH, HGG, KG, ZPF}, which
 has led to numerous research inventions.
Since 1950's, much attention has been devoted to modeling
and analyzing
chemostat problems; see \cite{FHC05,HEL56,HHW77,SW} and references therein.

The dynamics of the process can be modeled in the general form by a system of
ordinary differential equations
\begin{equation}\label{ww1}
\begin{cases}
\begin{aligned}
\dfrac{dS(t)}{dt}=&\dfrac{S_0-f_0(S(t))}\theta-X(t)f_1\big(S(t),X(t)\big),\\
\dfrac{dX(t)}{dt}=&X(t)\left(f_2\big(S(t),X(t)\big)-k_d-\dfrac{1+R}
\theta\right),
\end{aligned}
\end{cases}
\end{equation}
where
$S(t), X(t)$ are  the substrate concentration and
the bacterial concentration, respectively;
$S_0$ and $f_0(S)$
are the input concentration and the decay rate of the substrate, respectively;
$\frac 1\theta$ is the dilution rate
(or equivalently,
$\theta$ is the mean residence time),
$k_d$ is the death rate of $X$ and
$R$ is the recycle ratio;
$f_1(S,X)$ is the consumption rate and
$f_2(S,X)$ is the growth rate of the bacteria.
The formulation in \eqref{ww1}
is much more general
than the existing models.
The readers can find works on the specialized forms of the models with specific forms of the functions $f_0,f_1,f_2$ in \cite{BW85,HHW77,Zha16, ZY19,Zhi04} and references therein.

To better reflect the reality,
  effort has been devoted to stochastic systems to take into account the effect of
environmental perturbations \cite{IW, WJ}.
A fundamental problem is the long-term behavior of the system.
However, the dynamic behaviors have not been fully understood  to the best our knowledge. The asymptotic features of the systems and
important information such as the {\it wash-out time}
have not been
fully understood to date.
 In contrast to the existing work, we develop new approaches
to carefully analyze the corresponding systems.

Considering the system in a fluctuating environment,
we may assume that the dynamics are perturbed by a white noise.
Then, we have a
 stochastic counterpart of \eqref{ww1},
\begin{equation}\label{ww2}
\begin{cases}
\begin{aligned}
dS(t)=&\left(\dfrac{S_0-f_0\big(S(t)\big)}\theta-X(t)f_1\big(S(t),X(t)\big)\right)dt+\sigma_1S(t)dW_1(t),\\
dX(t)=&X(t)\left(f_2\big(S(t),X(t)\big)-k_d-\dfrac{1+R}\theta\right)dt+\sigma_2X(t)dW_2(t),
\end{aligned}
\end{cases}
\end{equation}
where $W_1(t)$ and $W_2(t)$ are two independent real-valued Brownian motions.
However, it has been recognized that the formulation above
is not able to capture
some important features of the underlying process. More often than not, in addition to the Brownian type perturbations,
there are also abrupt changes in the environment
that cannot be described by continuous perturbations.
An effective way to model these
discrete event
perturbations
is to use a Markov chain with a finite state space.
Suppose that the
coefficients $f_0,f_1,f_2,k_d$
and the intensities of the white noises
depend on
 $\alpha(t)$, a random switching process having a finite state space. Then
we have a more general system
\begin{equation}\label{ww3}
\begin{cases}
\begin{aligned}
dS(t)=&\left(\dfrac{S_0-f_0(S(t),\alpha(t))}\theta-X(t)f_1\big(S(t), X(t),\alpha(t)\big)\right)dt+\sigma_1(\alpha(t))S(t)dW_1(t),\\
dX(t)=&X(t)\left(f_2\big(S(t), X(t),\alpha(t)\big)-\wdt k_d(\alpha(t)) \right)dt
+\sigma_2(\alpha(t))X(t)dW_2(t),
\end{aligned}
\end{cases}
\end{equation}
where $\wdt k_d(\alpha(t)):= k_d(\alpha(t))+\frac{1+R}\theta$,  $\alpha(t)$ is a Markov chain with state space $\M=\{1,\dots, m_0\}$ and generator
$Q=(q_{kl})_{m_0\times m_0}$, and $\alpha(t)$ is independent of the
Brownian motions so that
\begin{equation}\label{eq:tran}\begin{array}{ll}
&\disp \PP\{\alpha(t+\Delta)=j|\alpha(t)=i,
\alpha(s), s\leq t\}=q_{ij}\Delta+o(\Delta) \text{ if } i\ne j \
\hbox{ and }\\
&\disp \PP\{\alpha(t+\Delta)=i|\alpha(t)=i,
\alpha(s), s\leq t\}=1+q_{ii}\Delta+o(\Delta).\end{array}\end{equation}

The models to be considered in this paper belong to the class of stochastic Kolmogorov systems, which is a class of dynamic systems used extensively in ecological, biological, and environment modeling.
Two long standing and fundamental questions concerning Kolmogorov systems are: (i) Under what conditions do populations persist or go extinct? (ii) When do interacting species coexist? The answers to these questions are essential for guiding conservation efforts. Our current paper is one of them along this line. In the literature, much effort has been devoted to such systems. A common approach is to use Lyapunov functions, which gives only sufficient conditions that are nowhere near necessary, and are not sharp. In addition, there is no systematic way of finding the Lyapunov functions. The analysis in this work is based on a completely different approach (used in \cite{MB, SBA11, HN18}), which requires treating the systems by looking at the boundary.
While sharp results for a very general class of stochastic systems have been obtained in \cite{MB}, some of their conditions are not always satisfied for our model. The proofs for persistence and extinction therefore require some delicate treatment of the behavior of the process near infinity. In the second part of the paper, we also consider controlled systems to reach the goal of getting optimality under a long-run average performance measure.

The rest of the paper is organized as follows.
In Section \ref{sec:thr}, we prove the existence and uniqueness of positive solutions to \eqref{ww3} and \eqref{eq:tran}.
Then a complete characterization of
the asymptotic behavior of the system under consideration
is provided. We show that the long-term properties of the system
can be classified by using a real-valued parameter $\lambda$.
If $\lambda\leq 0$,
the bacteria
will die out;
if $\lambda>0$, the system has an invariant probability measure
 and the transition probability of the solution process converges to  the invariant measure.
One of the distinct contributions of this paper is that the critical case $\lambda=0$ is also considered.
Some numerical examples are given in Section \ref{sec:num}.
Section \ref{sec:control} is devoted to the study of the system under ergodic control.
Controlled diffusions with random switching can be considered, but the notation will be more complex. To highlight the main ergodic control ideas, we decide to use a simplified model without switching.
We obtain the Hamilton-Jacobi-Bellman
 (HJB) equation and prove the existence and uniqueness of the solution of the HJB equation
corresponding to the long-term time-average control problem.
Establishing
 the existence and uniqueness of the solution to the HJB equation
is
most difficult
because the
usual conditions for ergodic controlled diffusions
are not satisfied in our model.
Inspired by the work \cite{ABG},
we use
a vanishing discount argument to
 examine
the associated cost and value functions
of the corresponding discounted control problem
and then to take a limit when the discount factor tends to $0.$
 However, the results in the
  aforementioned book cannot be applied or adopted directly because the conditions in the reference are not satisfied in
    our
  setup. More details on this will be given in Section \ref{sec:control}.
Finally, Section \ref{sec:con} issues  some concluding remarks.
Although our main motivation comes from
chemostat models, the analysis and techniques used can be applied to many other nonlinear ecological and biological systems.

\section{Complete Characterization of Long-Time Behavior}
\label{sec:thr}
This section is devoted to asymptotic properties of the switching diffusion models.
We show that
the long-term properties of the system
can be completely classified by using a real parameter $\lambda$.
In particular, if $\lambda\leq 0$, the bacteria
will die out and
we refer to such case as the system being not permanent because it does not work
in long term and the efficiency goes to 0.
If $\lambda>0$, the system has an invariant probability measure
 and the transition probability of the solution process converges to
the invariant measure.
This is what we refer to as permanence (with the terminology carried over from the study in biological and ecological models).
Moreover, we obtain the rate of convergence.

One of the highlights here is that we derive sufficient and necessary conditions for permanence.
First, the process under consideration is jointly a Markov Feller process. By using the Lyapunov exponent of $X(t)$, we then establish the existence of the invariant measure of $(S(t),\alpha(t))$. Furthermore,  under suitable conditions, we obtain the
exponential error bounds of the difference of the transition function and that of the invariant measure in the total variation norm.

Throughout this paper, we use the lowercase letters $s$, $x$, and $i$ to denote the initial values of $S(t)$, $X(t)$, and $\alpha(t)$, respectively.
Note the distinction of $s$ and
 the input concentration of the substrate $S_0$. To simplify the  notation, let
$$\hat k_d=\max_{i\in\M}\{\wdt k_d(i)\}\;;\;\check k_d=\min_{i\in\M}\{\wdt k_d(i)\},$$
and
 $$\hat \sigma_k=\max_{i\in\M}\{\sigma_k(i)\}\;;\;\check\sigma_k=\min_{i\in\M}\{\sigma_k(i)\}, k=1,2.$$
We denote by
$\R_+=[0,\infty)$, $\R^\circ_+=(0,\infty)$,
$\R_+^2=[0,\infty) \times [0,\infty)$, and $\R^{2,\circ}_+=(0,\infty)\times (0,\infty)$.
The operator associated with the process $\big(S(t),X(t),\alpha(t)\big)$, solving \eqref{ww3} and \eqref{eq:tran}, is given by
\begin{equation}
\Lom V(\phi,i)=V_\phi(\phi, i)\widetilde f(\phi, i)+\frac12
\tr [ \widetilde g(\phi, i) \widetilde g^\top(\phi,i) V_{\phi\phi}(\phi,i) ]
+\sum_{j\in\M} q_{ij}V(\phi, j),
\end{equation}
where
$A^\top$ denotes the transpose of $A$,
$\phi =(s,x)$,
$V_\phi(\phi, i) $ and $V_{\phi\phi}(\phi, i)$ are the gradient and Hessian of $V(\cdot,i)$ with respect to $\phi$,
$\widetilde f$ and $\widetilde g$ are the drift and diffusion coefficients of \eqref{ww3}, respectively.
That is,
\bea \ad \widetilde f(\phi,i)=\left(\dfrac{S_0-f_0(s,i)}\theta-xf_1(s,x,i),
\;\;x\Big(f_2(s,x,i)-\wdt k_d(i)\Big)\right)^\top,\\
 \ad \widetilde g(\phi,i)=\diag\big(\sigma_1(i) s, \sigma_2(i) x\big) \in \rr^{2\times 2},\eea
where
$\diag(a,b)$ denotes the diagonal matrix with entries $a$ and $b$. Note that the particular structure of $\wdt g$ implies that  $\wdt g \wdt g^\top = \wdt g^2$.
In what follows, we write $V(\phi,i)$ and $V(s,x,i)$ interchangeably,  whichever is more convenient. We also  assume the following conditions throughout.

\begin{asp}\label{as:1.1} {\rm
Suppose that
\begin{itemize}
\item $W_1(t)$ and $W_2(t)$ are independent real-valued standard Brownian motions
     that are independent of the Markov chain $\alpha(t)$.
\item  $S_0$, $\theta$, $\wdt k_d(i)$, $\sigma_1(i)$, $\sigma_2(i)$ are positive constants for each $i\in\M$.
\item $f_0(s,i)$, $f_1(s, x, i)$, and $f_2(s, x, i)$ are locally Lipschitz; $f_0(0,i)=f_1(0,x,i)=0$ satisfying
$\disp\lim_{s\to\infty} f_0(s,i)=\infty$;
$0\leq f_2(s, x, i)\leq \kappa_0 \big(f_1(s, x, i)\wedge f_0(s,i)\big)$,
$f_1(s,x,i)\leq \kappa_0(1+s)$
for some $\kappa_0>1$. Moreover, for each $i\in\M$, $f_2(s, x,i)$ is uniformly continuous at $x=0$ (that is,
$\lim_{x\to0}\sup\{|f_2(s,x,i)-f_2(s,0,i)|\}=0$).
\item The Markov chain or its generator
$Q=(q_{ij})_{m_0\times m_0}$ is irreducible, that is for any $i, j\in\M$, there exist $i=i_0, i_1, \dots, i_n=j$ such that $q_{i_{k-1},i_{k}}>0$ for $k=1,\dots,n$.

\end{itemize}
}
\end{asp}
To start, we state following Theorem, which provides the preliminary results, namely, the existence and uniqueness of the global solution,
and the positivity of the solution. The proof is postponed to the Appendix.

\begin{thm}\label{thm2.1}
For any $(s, x, i)\in\R_+^2\times\M$, there exists a unique global solution
to the system \eqref{ww3} and \eqref{eq:tran} with initial value $(s,x,i)$.
The three-component process $\{(S(t), X(t), \alpha(t)), t\geq0\}$
is a Markov-Feller process. Moreover, we have
$\PP_{s,x,i}\{S(t)>0, t>0\}=1$ and
$\PP_{s,x,i}\{X(t)=0, t>0\}=1$ if $x=0$,
$\PP_{s,x,i}\{X(t)>0, t>0\}=1$ if $x>0$.
\end{thm}

To proceed, we examine the boundary by letting $X(t)=0$.
Let $\hat S(t)$ be the solution to \eqref{ww3}
with $X(t)=0$, i.e.,
\begin{equation}\label{ww5}
d\hat S(t)=\dfrac{S_0-f_0\big(\hat S(t),\alpha(t)\big)}{\theta}dt+\sigma_1(\alpha(t))\hat S(t) dW_1(t).
\end{equation}
Since the drift of \eqref{ww5} is negative when $\hat S(t)$ is sufficiently large, and is positive when $\hat S(t)$ is sufficiently small,
we can
see
$$\left[\hat\op \Big(s+\ln\frac{s}{s+1}\Big)\right](s,i)\leq -1\text{ if } s<\hat s \text{ or } s>\frac1{\hat s},$$
for some $\hat s>0$ being sufficiently small,
where $\hat\op$ is the operator associated with \eqref{ww5}.
Thus, the hybrid diffusion \eqref{ww5} is positive recurrent due to its nondegeneracy (see, e.g., \cite[Chapter 4]{YZ}). Then there exists a unique invariant measure $\pi$ of \eqref{ww5}.
Moreover,
\begin{equation}\label{e2.11}
\lim_{t\to\infty}\|\hat P(t, s, i,\cdot)-\pi(\cdot)\|_{TV}=0, (s,i)\in(0,\infty)\times\M,
\end{equation}
where
$\|\cdot\|_{TV}$ is the total variation norm of a measure
and $\hat P(t,s,i,\cdot)$ is the transition probability of $(\hat S(t),\alpha(t))$.
Since $\PP_{0,i}\{\hat S(t)>0,t>0\}=1$,
\eqref{e2.11} holds even when $s=0$.

When $X(t)$ is small, $S(t)$ can be approximated by $\hat S(t)$.
Due to It\^o's formula and the ergodicity of $(\hat S(t),\alpha(t))$, when $X(t)$ is small,
the long-term growth rate $\dfrac{\ln X(t)}t$  is approximated by the critical value
\begin{equation}\label{lambda1}
\lambda
:=\sum_{i\in\M}\int_{\R_+^\circ}\left( f_2(s,0,i)-\wdt k_d(i)-\frac{\sigma_2^2(i)}2\right)\pi(ds,i).
\end{equation}
As a result, the sign of $\lambda$ determines whether or not $X(t)$ converges to $0$. A
 heuristic argument for getting $\lambda$ can be found in \cite{DNDY}. The main results of this section is provided in the next theorem to follow.
First, we state an assumption. One of the conditions
 in Assumption \ref{asp3} will be used in the main theorem.

\begin{asp}\label{asp3} {\rm Denoting by $\pi_\alpha$ the invariant measure of $\alpha(t)$,
assume that one of the following conditions holds.
\begin{enumerate}
\item[(1)] For each $i\in\M$, $f_2(s,x,i)$ is non-decreasing in $s$ and non-increasing in $x$
    for any $(s,x)\in[0,\infty)\times(0,\infty)$;
\item[(2)] \begin{equation}\label{cd2.1-e1}
\sum_{i\in\M}\left(\limsup_{s\to\infty} f_2(s,0,i)-\wdt k_d(i)-\frac{\sigma_2^2(i)}2\right)\pi_\alpha(i)<0;
\end{equation}
\item[(3)]
\begin{equation}\label{cd2.2-e1}
\sum_{i\in\M}\left(\liminf_{s\to\infty} f_2(s,0,i)-\wdt k_d(i)-\frac{\sigma_2^2(i)}2\right)\pi_\alpha(i)>0.
\end{equation}
\end{enumerate}
}\end{asp}

\begin{thm}\label{thm2.2} The following claims hold.
\begin{itemize}
\item If $\lambda<0$, we have
$\lim\limits_{t\to\infty}\dfrac{\ln X(t)}t=\lambda$ a.s. and  the distribution of $\big(S(t),\alpha(t)\big)$ converges weakly to
 the unique invariant probability measure $\pi$
 if either condition {\rm(1)} or {\rm(2)} of Assumption {\rm\ref{asp3}} holds.

 \item If $\lambda=0$ and condition {\rm(1)} of Assumption {\rm\ref{asp3}} holds  then
 \begin{equation}\label{2-e0}
 \lim_{T\to\infty} \dfrac1T\E_{s,x,i}\int_0^T S(t)dt=S_0\,\text{ and }\, \lim_{T\to\infty} \dfrac1T \E_{s,x,i}\int_0^T X(t)dt=0
 \end{equation}
 under additional conditions that for any $i\in\M$, $\frac{\partial f_0(s,i)}{\partial s}$ and $\frac{\partial f_2(s,0,i)}{\partial s}$ exist and are positive, continuous, and $\frac{\partial f_0(s,i)}{\partial s}$ is bounded below by a positive constant.
 \item If $\lambda>0$, then there exists an invariant probability measure $\mu^*$ on $\R^{2,\circ}_+\times\M$ and $\lim_{t\to\infty} \|P(t, s, i,\cdot)-\mu^*(\cdot)\|_{TV}=0$.
We assume further that  condition {\rm(3)} of Assumption {\rm\ref{asp3}} holds,
and $\liminf_{s\to\infty}\frac{f_0(s,i)}{s^q}>0, i\in\M$ for some $q\in(0,1]$.
Then,
\begin{itemize}
\item[\rm{(i)}]
in case $q<1$,
\bea \ad
\lim_{t\to\infty} t^{\beta-1}
 \|P(t, s,x, i,\cdot)-\mu^*(\cdot)\|_{TV}=0,\\
  \aad \hspace*{1.8in} \text{ for any }1\leq \beta<\frac{1}{1-q}, (s,x, i)\in\R_+\times\R^{\circ}_+\times\M;
 \eea
\item [\rm{(ii)}]
if $q=1$,
 there exists a $\wdt\gamma>0$  such that
 \begin{equation}\label{thm2.2-e0}
\lim_{t\to\infty} e^{\wdt\gamma t} \|P(t, s,x, i,\cdot)-\mu^*(\cdot)\|_{TV}=0, (s,x, i)\in\R_+\times\R^{\circ}_+\times\M,
 \end{equation}
 \end{itemize}
where $P(t, s,x, i,\cdot)$ is the transition probability of $\big(S(t),X(t),\alpha(t)\big)$.
 \end{itemize}
 \end{thm}

 \begin{rem}
 If $\lambda>0$, it is easy to see that condition (1) of Assumption \ref{asp3}
 implies the condition (3) in Assumption \ref{asp3}.
\blue{Moreover, most models in the literature consider the case
$f_0(s)=C_is$ while the coefficients $xf_1(s,x,i),xf_2(s,x)$ are linear functional response $(f_j(s,x,i)=c_{ij}x, j=1,2)$,  Holling type II response  $(f_j(s,x,i)=\frac{c_{ij}s}{m_{ij}+a_{ij}s}, j=1,2)$,
Holling type III response
$(f_j(s,x,i)=\frac{c_{ij}s^2}{(m_{ij1}+a_{ij1}s)(m_{ij2}+a_{ij2}s)}, j=1,2)$, and
 Beddington-DeAngelis functional response $(f_j(s,x,i)=\frac{c_{ij}s}{m_{ij}+a_{ij}s+b_{ij}x}, j=1,2)$
etc.;
see e.g., \cite{Zhi04,Zha16, ZY19} and the references therein. It is clear that these functions satisfy Assumption 2.1 and part (1) of Assumption 2.2.}
 \end{rem}

\begin{proof}[Proof for the case $\lambda=0$.]
 We argue by contradiction.
Suppose $\big(S(t),X(t),\alpha(t)\big)$ has an invariant probability measure $\mu^*$ on $\R^{2,\circ}_+\times\M$.
Then, it follows from the ergodicity that for $\mu^*$-almost every initial value $(s,x,i)$,
\begin{equation}\label{2-e1}
\lim_{T\to\infty}\frac1T\int_0^T h\big(S(t), X(t), \alpha(t)\big)dt
=\sum_{i'\in\M}\int_{\R^{2,\circ}_+}
h(s',x',i')\mu^*(ds',dx',i'),
\end{equation}
for any measurable function $h$ that is $\mu^*$-integrable.
Since the transition probability density is continuous and positive,
the invariant measure $\mu^*$ is unique and equivalent to $\wdt m\times \pi_\alpha$ where
$\wdt m$ is the Lebesgue measure on $\R^{2,\circ}_+$.
As a result, if $\mu^*(A)=1$, then $\PP_{s,x,i}\big((S(t), X(t),\alpha(t))\in A)\big)=1$ for every $(s,x,i)\in\R^{2,\circ}_+\times\M, t>0$,
which implies that \eqref{2-e1} holds for any initial value $(s,x,i)\in\R^{2,\circ}_+\times\M$.
By the comparison theorem (see e.g., \cite{GM94}), we have $\hat S(t)\geq S(t)$ with probability 1 given that $\hat S(0)=S(0)$. Hence, it follows \eqref{ww3}, \eqref{ww5} and positivity of $S(t)$, $\hat S(t)$ and additional assumption on $f_0(\cdot,\cdot)$ (non-decreasing in $s$) that
\begin{equation}\label{2-ebd}
\limsup_{T\to\infty} \dfrac 1T\E_{s,x,i} \int_0^T f_0(S(t),\alpha(t))dt\leq \limsup_{T\to\infty} \dfrac 1T\E_{s,i} \int_0^T f_0(\hat S(t),\alpha(t))dt\leq S_0.
\end{equation}
On the other hand, since $\dfrac{\partial f_0(s,i)}{\partial s}$ is bounded below by a positive constant and $f_0(0,i)=0$, there exists a $\bar g_1>0$ such that
\begin{equation}\label{f0'}
\frac{f_0(s,i)}{\theta}\geq \bar g_1s, \
 \forall s\geq 0, i\in\M.
\end{equation}
Let $0<\bar p<\frac{\bar g_1}{2\hat\sigma_1^2}$ and $0< \bar g_2<\frac{(1+\bar p)\bar g_1}{4}$. By It\^o's formula, we obtain that
$$
\begin{aligned}
\E_{s,i} &e^{\bar g_2T}\hat S^{1+\bar p}(T)\\
=&s^{1+\bar p}+\E_{s,i}\int_0^T e^{\bar g_2t}(1+\bar p)\hat S^{\bar p}(t)\left(\dfrac{S_0-f_0\big(\hat S(t),\alpha(t)\big)}{\theta}+\dfrac{\bar p\sigma_1^2(\alpha(t))}{2}\hat S(t)+\dfrac{\bar g_2}{1+\bar p}\hat S(t)\right)dt\\
\leq &s^{1+\bar p}+\E_{s,i}\int_0^T e^{\bar g_2t}(1+\bar p)\hat S^{\bar p}(t)\left(\dfrac{S_0}{\theta}-\dfrac{\bar g_1}2\hat S(t)\right)dt\\
\leq & s^{1+\bar p} +\dfrac{\bar g_3e^{\bar g_2T}}{\bar g_2},\;\forall T>0,\text{where }\bar g_3:=\dfrac{(1+\bar p)S_0}{\theta}\left(\dfrac{2S_0}{\theta\bar g_1}\right)^{\bar p}.
\end{aligned}
$$
As a consequence, we have
\begin{equation}\label{Shat-bd}
\limsup_{T\to\infty} \E_{s,x,i}S^{1+\bar p}(T)\leq\limsup_{T\to\infty}\E_{s,i}\hat S^{1+\bar p}(T)\leq \dfrac{\bar g_3}{\bar g_2}.
\end{equation}
Moreover,
we  also obtain that
\begin{equation}\label{lm2.1-e1}
\E_{s,x,i}\big(S(T)+X(T)\big)^{1+\bar p}\text{ is uniformly bounded in }T.
\end{equation}
Using It\^o's formula again, we have
\begin{equation*}
\begin{aligned}
\dfrac{\E_{s,i}\hat S^{1+\bar p}(T)}T=\dfrac{s^{1+\bar p}}{T}&+\dfrac{(1+\bar p)S_0}{\theta T}\E_{s,i}\int_0^T \hat S^{\bar p}(t)dt-\dfrac{1+\bar p}{\theta T}\E_{s,i}\int_0^T\hat S^{\bar p}(t)f_0\big(\hat S(t),\alpha(t)\big)dt\\
&+\dfrac{\bar p(1+\bar p)}{2T}\E_{s,i}\int_0^T \sigma_1^2(\alpha(t))\hat S^{1+\bar p}(t)dt,
\end{aligned}
\end{equation*}
which together with \eqref{Shat-bd} implies that
\begin{equation}\label{Sf0}
\limsup_{T\to\infty}\dfrac 1T\E_{s,i}\int_0^T\hat S^{\bar p}(t)f_0\big(\hat S(t),\alpha(t)\big)dt\leq \bar g_4,\;\text{for some }\bar g_4<\infty.
\end{equation}
We have from \eqref{ww5} and \eqref{Shat-bd} that
\begin{equation*}
\lim_{T\to\infty}\dfrac 1T\E_{s,i}\int_0^T \left(\dfrac{S_0}{\theta}-\dfrac{f_0\big(\hat S(t),\alpha(t)\big)}\theta\right)dt=\lim_{T\to\infty}\dfrac{\E_{s,i} \hat S(T)-s}{T}=0.
\end{equation*}
Hence
\begin{equation}\label{f0Shat}
\lim_{T\to\infty}\dfrac 1T\E_{s,i}\int_0^Tf_0\big(\hat S(t),\alpha(t)\big)dt=S_0.
\end{equation}

From \eqref{ww3}, with standard arguments, we have
$$\E_{s,x,i}\int_0^T X(t)f_1(S(t), X(t),\alpha(t))dt\leq S_0T+s,\; T\geq 0.$$
Moreover, due to the uniform boundedness \eqref{lm2.1-e1}, the linear growth bound of $f_0$ and the almost sure convergence \eqref{2-e1}, it follows from the dominated convergence theorem that
\begin{equation}\label{2-e7'}\lim_{T\to\infty}\frac1T\E_{s,x,i}\int_0^T\left(\dfrac{S_0-f_0\big(S(t),\alpha(t)\big)}{\theta}\right)dt =\sum_{i'\in\M}\int_{\R^{2,\circ}_+}
\frac{S_0-f_0(s',x',i')}\theta\mu^*(ds',dx',i')=:\hat g
\end{equation}
From \eqref{ww3} and \eqref{lm2.1-e1}, we deduces that $\lim_{T\to\infty}\frac1T\E_{s,x,i}\int_0^T X(t)f_1(S(t), X(t),\alpha(t))dt$ exists (because the two others in \eqref{2-e6} exist) and
\begin{equation}\label{2-e6}
\barray
\ad
\lim_{T\to\infty}\frac1T\E_{s,x,i}\int_0^T\left(\dfrac{S_0-f_0\big(S(t),\alpha(t)\big)}{\theta}-X(t)f_1\big(S(t),X(t),\alpha(t)\big)\right) dt\\
\aad \ =\lim_{T\to\infty}\frac{\E_{s,x,i}  S(T)-s}T= 0,\earray
\end{equation}
where, due to Fatou's lemma,
\bea
\ad
\lim_{T\to\infty}\frac1T\E_{s,x, i}\int_0^TX(t)f_1\big(S(t),X(t),\alpha(t)\big)dt \\
\aad \
=\hat g\geq \sum_{i'\in\M}\int_{\R^{2,\circ}_+}
x'f_1(s',x',i')\mu^*(ds',dx',i')>0.
\eea
We have from \eqref{f0Shat} and \eqref{2-e7'} that
\begin{equation}\label{2-e7}
\lim_{T\to\infty}\frac1T\E_{s,x, i}\int_0^T\left(f_0\big(\hat S(t),\alpha(t)\big)-f_0\big(S(t),\alpha(t)\big)\right) dt=\theta\hat g,\, s,x\geq0, i\in\M.
\end{equation}
Now, let $H>1$ be a sufficiently large constant satisfying that
$\frac{\bar g_4}{H^{\bar p}}\leq\frac{\theta\hat g}2. $
Combining with \eqref{Sf0}, we obtain
\begin{equation}\label{2-e8}
\begin{aligned}
\limsup_{T\to\infty}\dfrac 1T\E_{s,x, i}&\int_0^T\1_{\{\hat S(t)>H\}}\left(f_0\big(\hat S(t),\alpha(t)\big)-f_0\big(S(t),\alpha(t)\big)\right)dt\\
\leq& \limsup_{T\to\infty}\dfrac 1T\E_{s,x, i}\int_0^T\dfrac{\hat S^{\bar p}(t)}{H^{\bar p}}f_0\big(\hat S(t),\alpha(t)\big)dt\\
\leq& \dfrac {\bar g_4}{H^{\bar p}}
\leq\dfrac{\theta\hat g}2.
\end{aligned}
\end{equation}
A consequence of \eqref{2-e7} and \eqref{2-e8} is that
\begin{equation}\label{e2-11111}
\liminf_{T\to\infty}\frac1T\E_{s,x, i}\int_0^T\1_{\{\hat S(t)\leq H\}}\left(f_0\big(\hat S(t),\alpha(t)\big)-f_0\big(S(t),\alpha(t)\big)\right) dt\geq \frac{\theta\hat g}2.
\end{equation}
Since $\frac{\partial f_0(s,0,i)}{\partial s}$ is continuous and positive, there exists a $\bar g_5^H<\infty$ such that
$$f_0(s_1,i)-f_0(s_2,i)\leq \bar g_5^H(s_1-s_2),\;\forall 0\leq s_2\leq s_1\leq H, i\in\M.$$
Thus, \eqref{e2-11111} implies that
\begin{equation}\label{e2-11112}
\liminf_{T\to\infty}\frac1T\E_{s,x, i}\int_0^T\1_{\{\hat S(t)\leq H\}}\left(\hat S(t)-S(t)\right) dt\geq \frac{\theta\hat g}{2\bar g_5^H}.
\end{equation}
Similarly, there exists a $\bar g_6^H>0$ such that
$$f_2(s_1,0,i)-f_2(s_2,0,i)\geq \bar g_6^H(s_1-s_2),\text{ for any } 0\leq s_2\leq s_1\leq H, i\in\M.$$
 Hence, combining with $f_2(s,x,i)$ being non-increasing in $x$, we obtain
\begin{equation}\label{e2-11113}
f_2(s_1,0,i)-f_2(s_2,x,i)\geq \bar g_6^H(s_1-s_2)\;\forall x\geq 0,0\leq s_2\leq s_1\leq H, i\in\M.
\end{equation}
As a consequence of \eqref{e2-11112} and \eqref{e2-11113},
\begin{equation*}
\begin{aligned}
\lim_{T\to\infty}&\frac1T\E_{s,x, i}\int_0^T\left(f_2\big(\hat S(t),0,\alpha(t)\big)-f_2\big(S(t),X(t),\alpha(t)\big)\right) dt
\\
\geq&
\lim_{T\to\infty}\frac1T\E_{s,x, i}\int_0^T\1_{\{\hat S(t)\leq H\}}\left(f_2\big(\hat S(t),0,\alpha(t)\big)-f_2\big(S(t),X(t),\alpha(t)\big)\right) dt
\\
\geq&\bar g_6^H\liminf_{T\to\infty}\frac1T\E_{s,x, i}\int_0^T\1_{\{\hat S(t)\leq H\}}\left(\hat S(t)-S(t)\right) dt\\
\geq&\frac{\bar g_6^H\theta\hat g}{2\bar g_5^H}.
\end{aligned}
\end{equation*}
Therefore, we obtain that
\begin{equation}\label{04-23-1}
\begin{aligned}
\sum_{i'\in\M}&\int_{\R_+^{2,\circ}}\left( f_2(s',x',i')-\wdt k_d(i)-\frac{\sigma_2^2(i)}2\right)\mu^*(ds',dx', i')\\
=&\lim_{T\to\infty}\frac1T\E_{s,x, i}\int_0^T\left(f_2\big(S(t),X(t),\alpha(t)\big)-\wdt k_d(\alpha(t))-\dfrac {\sigma_2^2(\alpha(t))}{2}\right) dt\\
\leq &\lim_{T\to\infty}\frac1T\E_{s,x, i}\int_0^T\left(f_2\big(\hat S(t),0,\alpha(t)\big)-\wdt k_d(\alpha(t))-\dfrac {\sigma_2^2(\alpha(t))}{2}\right) dt-\frac{\bar g_6^H\theta\hat g}{2\bar g_5^H}\\
\leq& \lambda-\frac{\bar g_6^H\theta\hat g}{2\bar g_5^H}=-\frac{\bar g_6^H\theta\hat g}{2\bar g_5^H}.
\end{aligned}
\end{equation}
By \eqref{ww3}, It\^o's formula, the ergodicity, and \eqref{04-23-1}, we have
$$
\begin{aligned}
\lim_{T\to\infty}\dfrac{\ln X(T)}T
=&
\lim_{T\to\infty}\dfrac1T\left(X(0)+\int_0^T\sigma_2(\alpha(t))dW_2(t)\right)\\
&+\lim_{T\to\infty}\dfrac1T\int_0^T\left(f_2\big(S(t),X(t),\alpha(t)\big)-\wdt k_d(\alpha(t))-\dfrac {\sigma_2^2(\alpha(t))}{2}\right) dt\\
\leq& -\frac{\bar g_6^H\theta\hat g}{2\bar g_5^H}<0\ \text{ a.s.}
\end{aligned}
$$
As a result,
$$\PP_{s,x,i}\left\{\lim_{T\to\infty} X(T)=0\right\}=1,$$
which contradicts the assumption that the process has an invariant probability measure on $\R^{2,\circ}_+\times\M$.
Moreover,
$\pi\times\delta^*$ is the unique invariant measure of $\big(S(t),\alpha(t), X(t)\big)$ on $\R_+\times\M\times\R_+$,
where $\delta^*$ is the Dirac measure with mass at $0$.
Consider the empirical measure
$$\Pi^{s,x,i}_t(\cdot)=\dfrac1t\int_0^t\PP_{s,x,i}\{(S(s),X(s), \alpha(s))\in\cdot\}ds.$$
In view of \eqref{lm2.1-e1},
the family $\{\Pi^{s,x,i}_t(\cdot), t\geq 0\}$ is tight for each $(s,x,i)\in\R^2_+\times\M$.
It is well-known (see e.g., \cite{HN18} and \cite{SBA11}) that any weak limit of $\Pi^{s,x,i}_t(\cdot)$ as $t\to\infty$
is an invariant probability measure of $\big(S(t),\alpha(t), X(t)\big)$.
Since $\pi\times\delta^*$ is the unique invariant probability measure and we have the boundedness of $\E_{s,x,i}\big(S(t)+X(t)\big)^{1+\bar p}$ in \eqref{lm2.1-e1},
we can easily
obtain \eqref{2-e0}.
\end{proof}

\begin{proof}[Proof for the case $\lambda>0$]
	\blue{Proof of the convergence in total variation of transition probability to an invariant measure is
straightforward.
Applying \cite[Theorem 4.4]{MB},
with $W(s,x,i)=s+x$ and $V(s, x, i)=\ln\frac{x}{x+1}$
we obtain the persistence of $X(t)$ and the existence of an invariant probability measure of $(S(t), X(t),\alpha(t))$ on $\R^{2,\circ}_+\times\M$.
Because of the nondegeneracy of the diffusion, which implies irreducibility and strong Feller property of the skeleton process $\{(S(nt_0)$, $X(nt_0)$, $\alpha(nt_0)$, $n\in\N\}$
for any $t_0>0$,
we can obtain the convergence in total variation of transition probability of the process $(S(t), X(t),\alpha(t))$ to its  invariant probability measure on $\R^{2,\circ}_+\times\M$;
see \cite[Theorem 6.1]{MT2} or \cite{ZY}.
}

Now, we consider the rate of convergence where conditions in \cite[Theorem 4.12 or Proposition 4.13]{MB} are not easy to verify for our model. We assume that condition (3) of Assumption \ref{asp3} holds.
 In view of \eqref{cd2.2-e1}, there exists an $\hat H>0$ such that
\begin{equation}\label{lambda1b}
\sum_{i\in\M}\left( \inf_{s\geq\hat H}\{f_2(s,0,i)\}-\wdt k_d(i)-\frac{\sigma_2^2(i)}2\right)\pi_\alpha(i)>0.
\end{equation}
Therefore, by the uniform continuity at $x=0$ of $f_2(s,x,i)$, there exists an $\eps_1>0$ such that
\begin{equation}\label{lambda1b}
4\wdt\lambda:=\sum_{i\in\M}\hat h_i\pi_\alpha(i)>0, \text{ where } \hat h_i:= \inf_{(s,x)\in[\hat H,\infty)\times[0,\eps_1]}\{f_2(s,x,i)\}-\wdt k_d(i)-\frac{\sigma_2^2(i)}2.
\end{equation}
Since $\sum_{i\in\M}(4\wdt\lambda-\hat h_i)\pi_\alpha(i)=0$, an application of the Fredholm alternative implies the existence of $\gamma_i>0$ with $i\in\M$ such that
$$
\sum_{j\in\M} q_{ij}\gamma_j=4\wdt\lambda-\hat h_i,\;\forall i\in\M.
$$
Let $\rho_1$ be sufficiently small satisfying
$$\rho_1\gamma_i(-4\wdt\lambda+\hat h_i)<\wdt\lambda(1-\rho_1\gamma_i),\;\rho_1\sigma_2^2(i)<4\wdt\lambda\text{ and } \rho_1\gamma_i<1,\forall i\in\M,$$
and define $V_3(s,x,i)=(1-\rho_1\gamma_i)x^{-\rho_1}$. By directed calculations, we obtain that
\begin{equation}\label{e-LV3}
\begin{aligned}
\Lom V_3(s,x,i)=&\rho_1V_3(s,x,i)\left(-f_2(s,x,i)+\wdt k_d(i)+\dfrac{\rho_1+1}2\sigma_2^2(i)\right)\\
&+\rho_1(-4\wdt\lambda+\hat h_i)V_3(s,x,i)+\rho_1^2\gamma_i(-4\wdt\lambda+\hat h_i)x^{-\rho_1}\\
\leq & -2\rho_1\wdt\lambda V_3(s,x,i)+\rho_1^2\gamma_i(-4\wdt\lambda+\hat h_i)x^{-\rho_1}\\
\leq &-\rho_1\wdt \lambda V_3(s,x,i)\;\forall i\in\M, s\geq \hat H, x\leq \eps_1.
\end{aligned}
\end{equation}
Since $\left[\op \frac1x\right](s,x,i)\leq \frac{\wdt k_d(i)+\sigma_2^2(i)}x\leq \frac{\hat k_d+\hat \sigma_2^2}x$, it is easily verified that
\begin{equation}\label{thm2.2-e1a}
\E_{s,x,i} X^{-1}(t)\leq e^{(\hat k_d+\hat \sigma_2^2) t} x^{-1}, t\geq0, (s,x,i)\in\R^{2,\circ}_+\times\M.
\end{equation}
On the other hand, by \eqref{lambda1}, we obtain that there is a $T_1>1$ satisfying
\begin{equation}\label{thm2.2-e4a}
-\ln(1-\rho_1\gamma_i)<\frac{\rho_1\lambda T_1}{4},\forall  i\in\M,
\end{equation}
such that
\begin{equation}\label{e3.0}
\dfrac1T\E_{s, 0, i}\int_0^T\left(f_2\big(S(t), 0, \alpha(t)\big)-\wdt k_d(\alpha(t))-\frac{\sigma_2^2(\alpha(t))}2\right)dt>\frac{3\lambda}4, i\in\M, T\geq T_1, s\leq \hat H.
\end{equation}
Let $n_*\in\N$ and $n_*>\frac{\hat k_d+\hat \sigma_2^2}{\wdt\lambda}+1$.
Inspired by the use of the log-Laplace transform in \cite{MB,BL16},
we can follow the proof of \cite[Proposition 4.1]{HN18}
to obtain the existence of some $\rho_2\in(0,\rho_1)$ and $\eps_2\in(0,\eps_1)$ satisfying
\begin{equation}\label{thm2.2-e2a}
 \E_{s,x,i} X^{-\rho_2}(t)\leq e^{-\frac{\rho_2\lambda t}2}x^{-\rho_2}\text{ for } i\in\M, t\in[T_1, n_*T_1], s\leq \hat H, x<\eps_2.
 \end{equation}
With such $\rho_2$,
we have from \eqref{thm2.2-e2a}  and \eqref{thm2.2-e4a} that
\begin{equation}\label{thm2.2-e3a}
\begin{aligned}
\E_{s,x,i} &V_3^{\frac{\rho_2}{\rho_1}}\big(S(t), X(t), \alpha(t)\big)\\
\leq& \E_{s,x,i} X^{-\rho_2}(t)\leq e^{-\frac{\rho_2\lambda t}2}x^{-\rho_2}=e^{-\frac{\rho_2\lambda t}2} V_3^{\frac{\rho_2}{\rho_1}}(s,x,i)(1-\rho_1\gamma_i)^{-\frac{\rho_2}{\rho_1}}\\
\leq &e^{-\frac{\rho_2\lambda t}2} V_3^{\frac{\rho_2}{\rho_1}}(s,x,i) e^{\frac{\rho_2\lambda t}4}\leq e^{-\frac{\rho_2\lambda t}4} V_3^{\frac{\rho_2}{\rho_1}}(s,x,i),\forall i\in\M, t\in[T_1, n_*T_1], s\leq \hat H, x<\eps_2.
\end{aligned}
\end{equation}
Since $\rho_2<\rho_1$,
we can  apply It\^o's formula to obtain from \eqref{e-LV3}
that
\begin{equation}\label{thm2.2-e5a}
\op V_3^{\frac{\rho_2}{\rho_1}}(s,x,i)\leq -\rho_2\wdt\lambda V_3^{\frac{\rho_2}{\rho_1}}(s,x,i)\;\forall i\in\M, s\geq \hat H, x\leq \eps_2.
\end{equation}
Estimates
\eqref{thm2.2-e3a} and \eqref{thm2.2-e5a} allow us to estimate
$\E_{s,x,i} V_3^{\frac{\rho_2}{\rho_1}}\big(S(t), X(t), \alpha(t)\big)$ through $V_3^{\frac{\rho_2}{\rho_1}}(s,x,i)$
when $x$ is sufficiently small and $s\geq 0$.
Using \eqref{thm2.2-e3a}, \eqref{thm2.2-e5a}, and \eqref{thm2.2-e1a},
we can follow the proof of \cite[Theorem 4.1]{HN18} to show the existence of $q_*\in(0,1)$ and $K_*>0$ such that
\begin{equation}\label{v3-e3}
\E_{s,x,i} V_3^{\frac{\rho_2}{\rho_1}}\big(S(n_*T_1), X(n_*T_1), \alpha(n_*T_1)\big)
\leq  q^* V_3^{\frac{\rho_2}{\rho_1}}(s,x,i)+K_*.
\end{equation}
Moreover,
if $\liminf_{s\to\infty} \frac{f_0(s,i)}{s^q}>0, i\in\M$ for some $q\in(0,1]$,
we have
\begin{equation}\label{v3-e4}
\left[\op (\kappa_0 s+x)\right](s,x,i)\leq \frac{S_0\kappa_0}{\theta}- \frac{\kappa_0 f_0(s,i)}\theta-\wdt k_d(i) x
\leq \wdt c_1-\wdt c_2(\kappa_0 s + x)^q,
\end{equation}
for some $\wdt c_1,\wdt c_2>0$.

Having \eqref{v3-e3} and \eqref{v3-e4},
we can apply \cite[Theorem 3.6]{JR02} and \cite[Proposition 22]{FR05} with slight
modifications of their proofs to show that
\begin{equation}\label{v3-e5}
\E_{s,x,i} \sum_{k=1}^{\tau^*(\wdt {\mathcal H})}(k+1)^{\beta-1}\leq c_\beta(\wdt {\mathcal H})\left( V_3^{\frac{\rho_2}{\rho_1}}(s,x,i)+\kappa_0 s+ x+1\right),
\end{equation}
where $1\leq \beta<\frac{1}{1-q}$,  $\wdt {\mathcal H}$ is a compact set in $[0,\infty)\times(0,\infty)$,
$c_\beta(\wdt {\mathcal H})$ is some positive constant,
and $$\tau^*(\wdt {\mathcal H})=\inf\{k\in\N:  (S(kn_*T_1), X(kn_*T_1))\in\wdt {\mathcal H}\}.$$
Since the drift term of $S(t)$ is positive when $S(t)=0$,
a standard argument shows that
the compact set $\wdt {\mathcal H}$ of $[0,\infty)\times(0,\infty)$  is petite (see e.g., \cite{DNDY}).
Using this fact and \eqref{v3-e5}, we derive from \cite[Theorem 2.1]{TT94} that
$$
\lim_{k\to\infty} (k+1)^{\beta-1}
 \|P(kn_* T_1, s, x,i,\cdot)-\mu^*(\cdot)\|_{TV}=0.
 $$
An application of some standard arguments
(see e.g., \cite{DNDY})
then shows that
$$
\lim_{t\to\infty} t^{\beta-1}
 \|P(t, s, x,i,\cdot)-\mu^*(\cdot)\|_{TV}=0,\text{ for any } 1\leq \beta<\frac{1}{1-q}.
 $$

If $q=1$, we obtain from \eqref{v3-e3} and \eqref{v3-e4} that
\begin{align*}
\E_{s,x,i} &\left[V_3^{\frac{\rho_2}{\rho_1}}\big(S(n_*T_1), X(n_*T_1), \alpha(n_*T_1)\big)+\kappa_0 S(n_*T_1)+X(n_*T_1)\right]
\\&\leq  \wdt q^* \left(V_3^{\frac{\rho_2}{\rho_1}}(s,x,i)+\kappa_0s+x\right)+\wdt K_*,
\end{align*}
for some $\wdt q^*\in(0,1)$ and $\wdt K_*>0$.
Then \cite[Theorem 2.1]{TT94} implies that
$$
\lim_{k\to\infty} e^{\wdt\gamma k}
 \|P(kn_* T_1, s,x, i,\cdot)-\mu^*(\cdot)\|_{TV}=0, \text{ for some } \wdt\gamma>0.
 $$
 Then,
we can  obtain the exponential rate of convergence.
\end{proof}

\begin{proof}[Proof for the case $\lambda<0$.]
If condition (1) of Assumption \ref{asp3} holds, the proof can be carried out
using comparison arguments by comparing $S(t)$ and $\hat S(t)$.
The details are omitted here since they are similar to  \cite[Theorem 2.1]{DNDY}.
[In \cite{DNDY}, we obtain that $\limsup_{t\to\infty}\frac{\ln X(t)}t\leq\lambda$, then simple arguments using the fact that any weak limit of random occupation measures is an invariant probability measure (see e.g., \cite{HN18}),
and
 we can obtain the
  convergence rate (Lyapunov exponent)
 $\lim_{t\to\infty}\frac{\ln X(t)}t=\lambda$.]

Now, we assume that condition (2) of Assumption \ref{asp3} holds.
Analogous to the proof in the case $\lambda>0$,
we have
$$\E_{s,x,i} V_5^{\frac{\rho_4}{\rho_3}}\big(S(n^*T_2), X(n^*T_2), \alpha(n^*T_2)\big)
\leq \hat q^* V_5^{\frac{\rho_4}{\rho_3}}(s,x,i)+\hat K_*;  i\in\M,s\geq 0,x\leq\eps_3,
$$
where
$V_5(s,x,i)=(1-\rho_3\gamma_i)x^{\rho_3}$
for suitable constants $\rho_3,\rho_4$, $n^*$, $T_2$, $\hat q^*\in(0,1), \hat K_*$, $\eps_3$.
Then we can mimic the proofs of \cite[Theorems 5.1 \&  5.2]{HN18} to obtain the desired result.
\end{proof}

\begin{rem} Before proceeding further, let us make the following remarks.
	\begin{itemize}
		\item
	\blue{The techniques of handling the critical case can be applied to
such
stochastic models in epidemiology
as SIR, SEIR, SEIS
			models, and
predator-prey models,
			when the system exhibits certain monotone properties
			so that the contradiction arguments in the proof for the case $\lambda=0$ can
be applied.	}
\item
\blue{Constructing Lyapunov function for switching diffusions
is practically more difficult than that for diffusions.
In non-critical cases, we showed how to construct the pair of Lyapunov functions  $(\wdt V, \wdt H)$ satisfying \cite[Proposition 4.13]{MB}
for systems involving Markovian switching
if the construction is possible but not obvious
(in our model, the requirement is
(2) or (3) of Assumption \ref{asp3} together with $\lim_{s\to\infty} \frac{f_0(s,i)}{s^q}>0, q<1$ being  satisfied).
We also showed that when it is practically impossible to find
  $(\wdt V, \wdt H)$ satisfying \cite[Proposition 4.13]{MB},
  we can obtain a sub-geometric convergence rate under certain conditions. Our techniques combined with the main approaches
in \cite{MB} can work for some stochastic Kolmogorov systems
under relaxed conditions.
  }
	\end{itemize}

\end{rem}

\section{Numerical Examples}\label{sec:num}

This section is devoted to
some numerical examples. Consider an
example of  system \eqref{ww3} as follows, which is
often
used in wastewater treatment (see e.g., \cite{NS,SW})
\begin{equation}\label{ww3-111}
\begin{cases}
\begin{aligned}
dS(t)=&\left(\dfrac{S_0-S(t)}\theta-\dfrac{k_m(\alpha(t))S(t)X(t)}{K_S+S(t)}\right)dt+\sigma_1(\alpha(t))S(t)dW_1(t),\\
dX(t)=&X(t)\left(\dfrac{k_m(\alpha(t))Y(\alpha(t))S(t)}{K_S+S(t)}-k_d(\alpha(t))-\dfrac{1+R}{\theta} \right)dt
+\sigma_2(\alpha(t))X(t)dW_2(t).
\end{aligned}
\end{cases}
\end{equation}
The following
 table provides parameter values for conventional activated sludge system using a completely mixed flow reactor extracted from \cite[p. 351]{NAC}.
\begin{center}
\begin{tabular}{ |l|l|l| }
 \hline
 Parameter & Typical range & Units \\
 \hline
 $k_m$: the growth constant of the bacteria & 2-10 & mg of cells $\times$ day \\
 $K_S$: the half-saturation constant & 25-100 & mg of cells $\times$ day/L \\
 $Y$=$\frac{\text{the growth rate of the bacteria}}{\text{the rate of substrate consumption}}$ & 0.4-0.8& dimensionless\\
 $k_d$: the death rate&0.025--0.075 &1/day\\
$\theta$: the hydraulic residence time& 3-5& day\\
  \hline
\end{tabular}
\end{center}

\

\begin{exam}\label{ex1}{\rm
Consider equation \eqref{ww3-111} with $\alpha(t)\in\M=\{1,2\}$ and parameters
$
S_0=15$,
$k_m(1)=9$,
$k_m(2)=6$,
$\theta=5$,
$R=0$,
$Y(1)=0.8$,
$Y(2)=0.6$,
$k_d(1)=0.06$,
$k_d(2)=0.08$,
$K_S=60$,
$\sigma_1(1)=0.1$, $\sigma_2(1)=0.2$,
$\sigma_1(2)=1$, $\sigma_2(2)=0.1$,
$q_{12}=0.2$, and $q_{21}=0.8$.
Using the strong law of large numbers,
$$\lambda=\lim_{T\to\infty}\frac1T\int_0^T
\left(\dfrac{k_m(\alpha(t))Y(\alpha(t))\hat S(t)}{K_S+\hat S(t)}-k_d(\alpha(u))\right).$$
We can approximate $\lambda$ through the occupation measure in a long period of time $[0,T]$.
In this example,
$\lambda\approx  0.915>0$.
Thus, the process $(S(t), X(t))$ has an invariant probability measure on $\R^{2,\circ}$.
Figure \ref{f1.1} displays
a sample path of $S(t), X(t)$.
The empirical approximation for the density function
is shown in Figures \ref{f1.2} and \ref{f1.3}.

   \begin{figure}[htp]
\centering
\includegraphics[width=7.5cm, height=7cm]{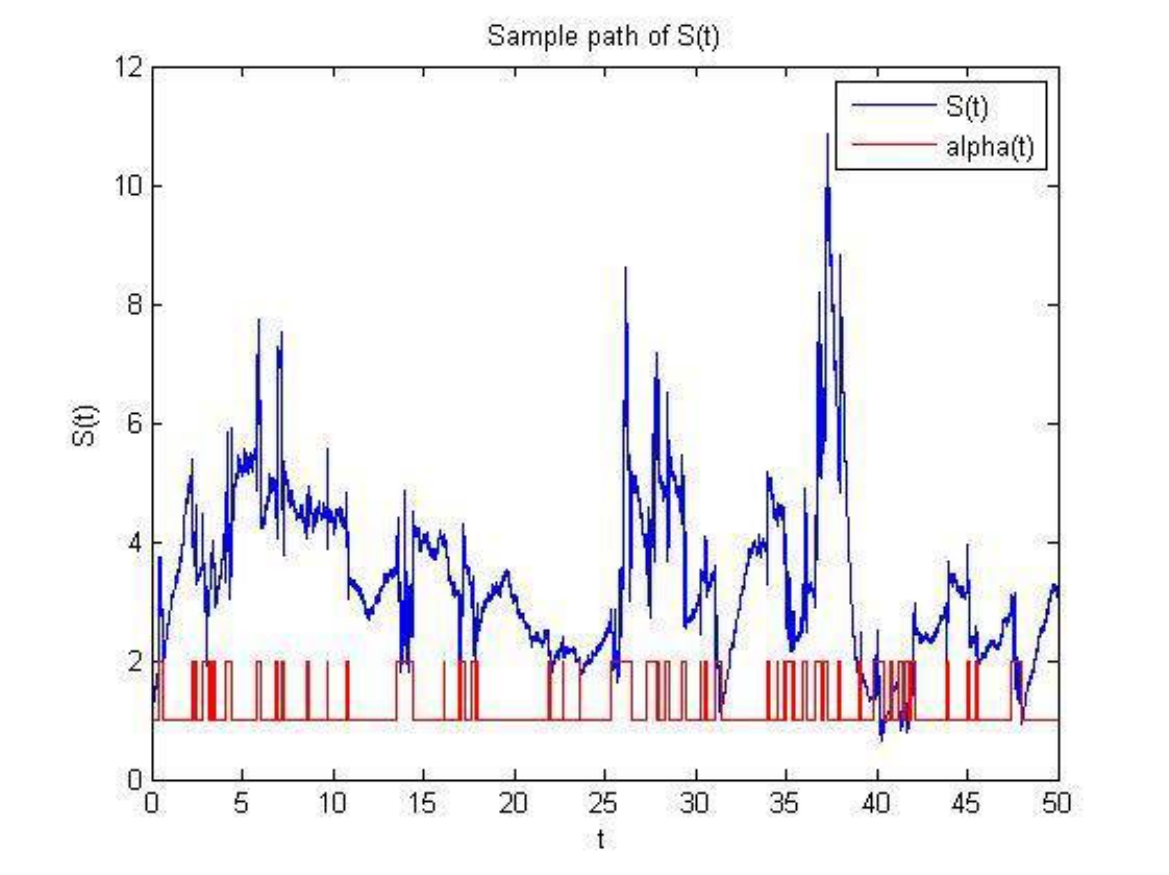}
\includegraphics[width=7.5cm, height=7cm]{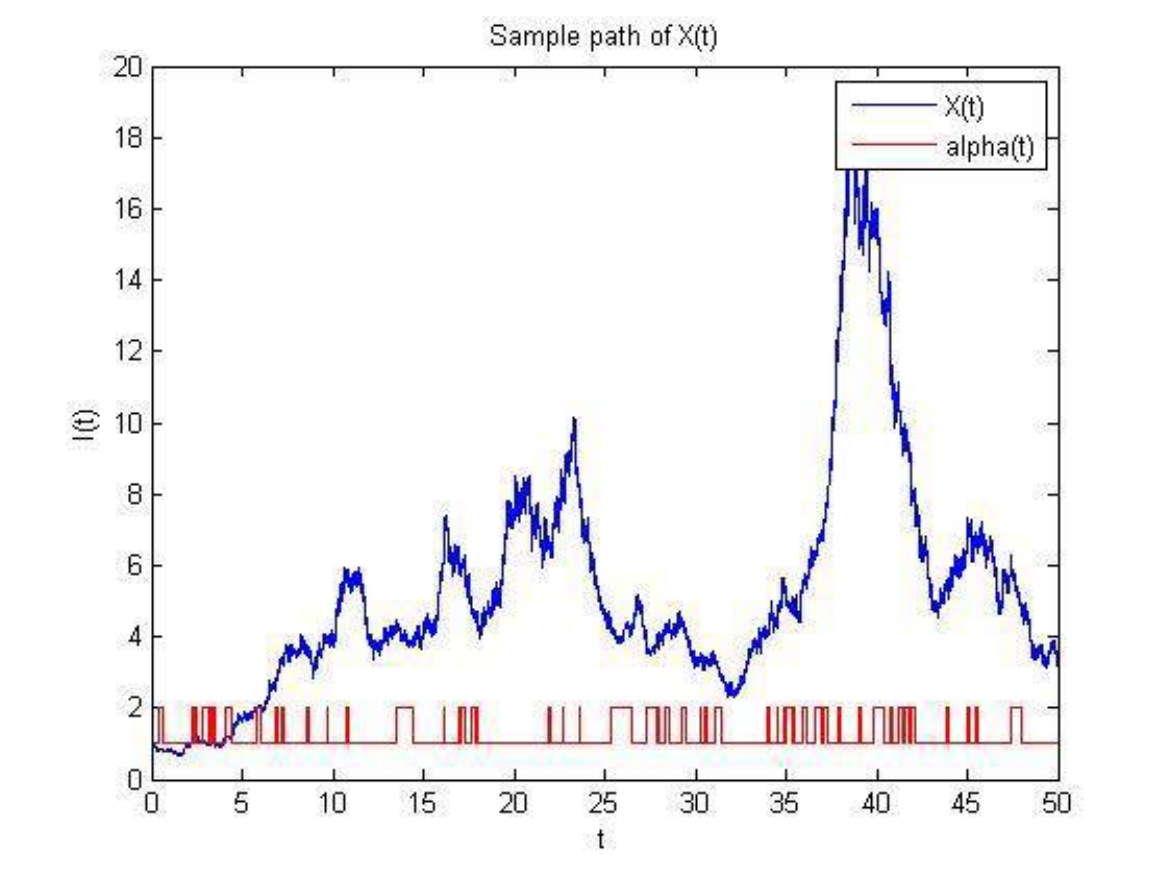}
\caption{Sample paths of $S(t)$ (in blue on the left) and $X(t)$ (in blue on the right) and $\alpha(t)$ (in red) in Example \ref{ex1}.}\label{f1.1}
\end{figure}

   \begin{figure}[htp]
\centering
\includegraphics[width=7.5cm, height=7cm]{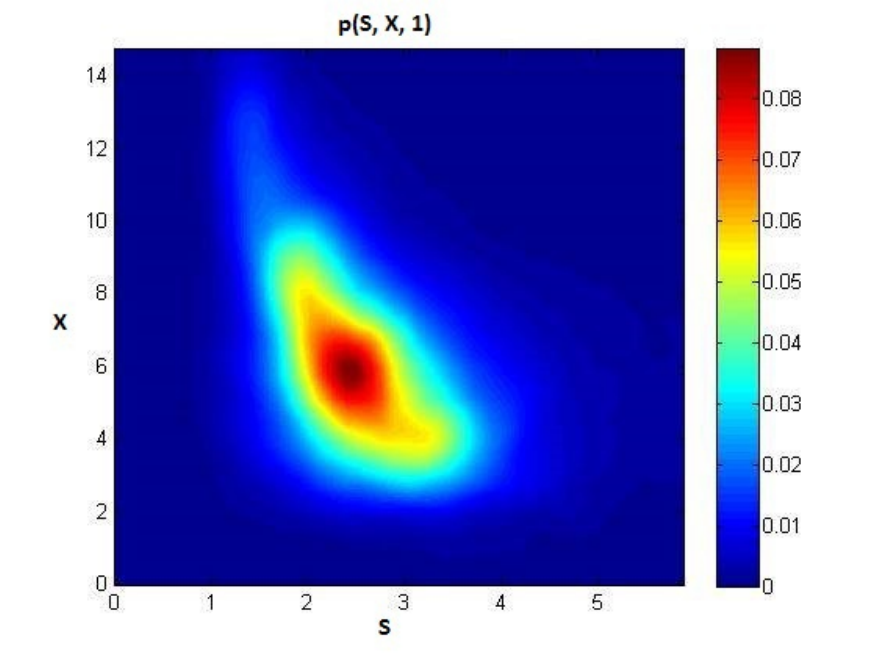}
\includegraphics[width=7.5cm, height=7cm]{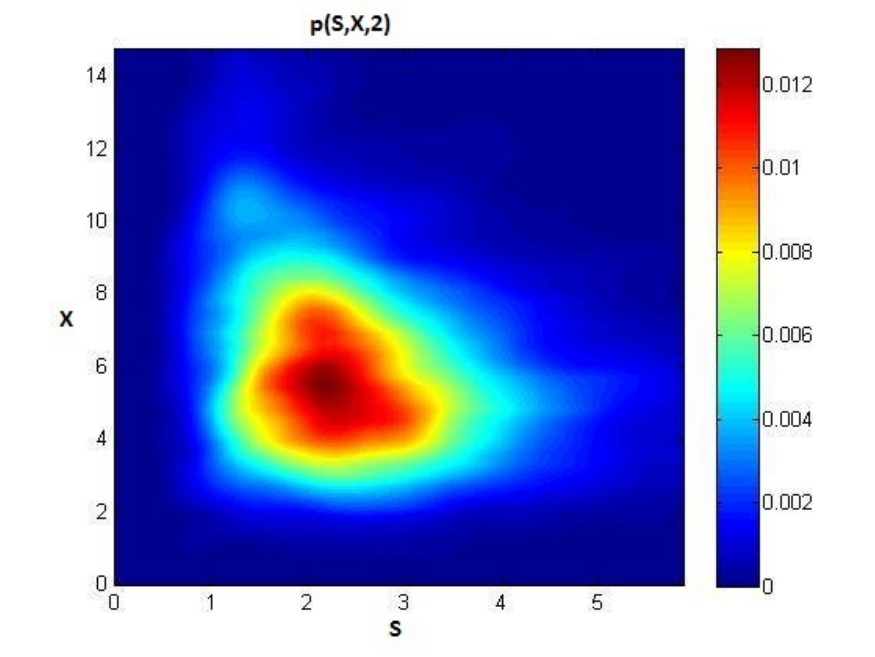}
\caption{
Densities $p(s,x,1)=\mu^*(ds,dx,1)$ (on the left) and $p(s,x,2)=\mu^*(ds,dx,2)$ (on the right) of the invariant probability measure $\mu^*$
 in Example \ref{ex1}. Different colors  represent different sizes of the density.
}\label{f1.2}
\end{figure}

   \begin{figure}[htp]
\centering
\includegraphics[width=7.5cm, height=7cm]{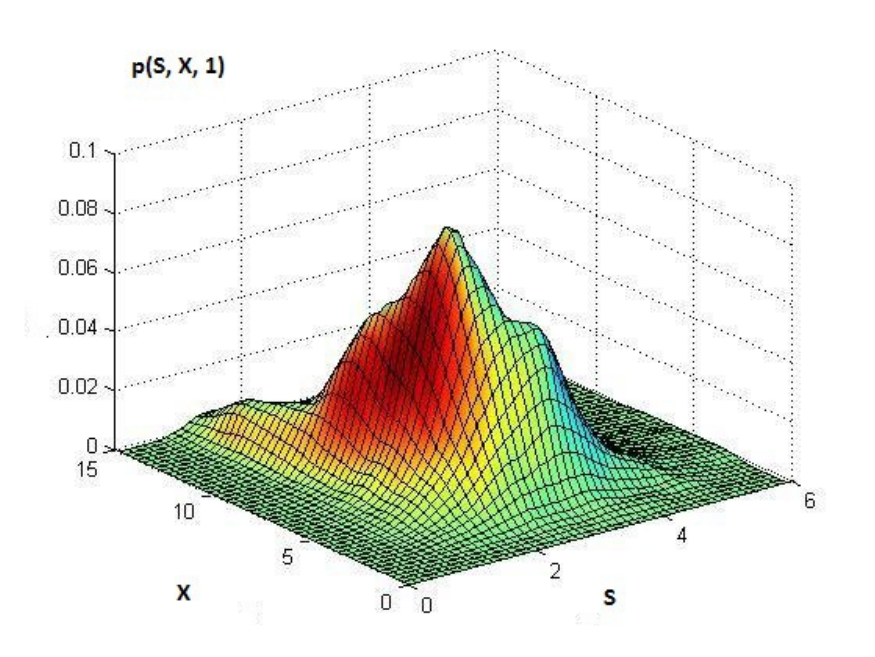}
\includegraphics[width=7.5cm, height=7cm]{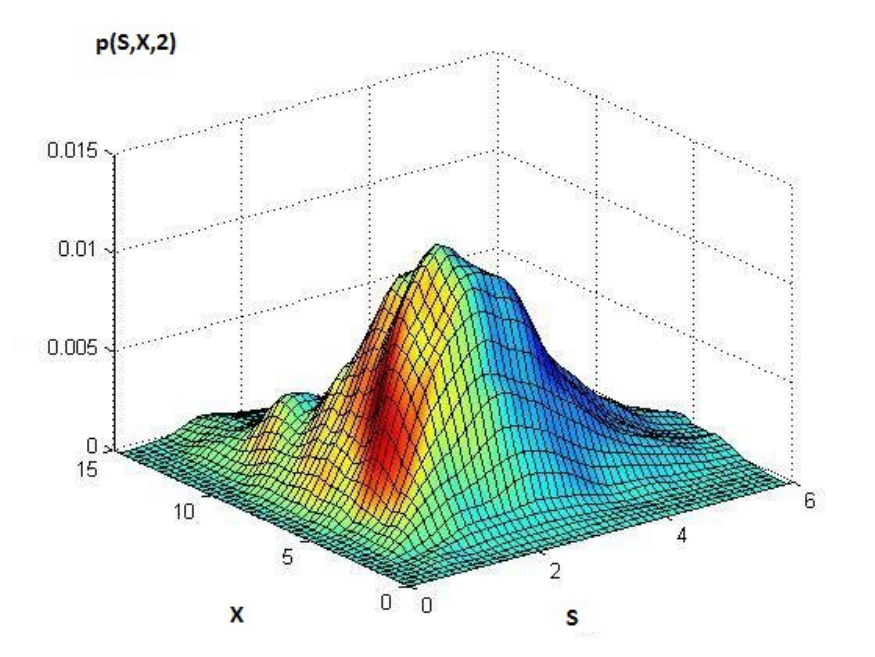}
\caption{
Graphs of densities $p(s,x,1)$ (on the left) and $p(s,x,2)$ (on the right) of the invariant probability measure $\mu^*$
 in Example \ref{ex1}.
}\label{f1.3}
\end{figure}

}\end{exam}

\begin{exam}\label{ex2}{\rm
Consider \eqref{ww3-111} without switching and parameters
$$
S_0=12,
k_m=8,
\theta=1,
\rho=0,
Y=0.6,
k_d=0.06,
K_S=60,
\sigma_1=0.2, \sigma_2=0.2.
$$
Direct
computation shows that $\lambda\approx -0.28<0$.
Thus, $X(t)$ will tend to $0$ as $t\to\infty$,
which is illustrated in Figure \ref{f2.1}.
   \begin{figure}[htp]
\centering
\includegraphics[width=7.5cm, height=7cm]{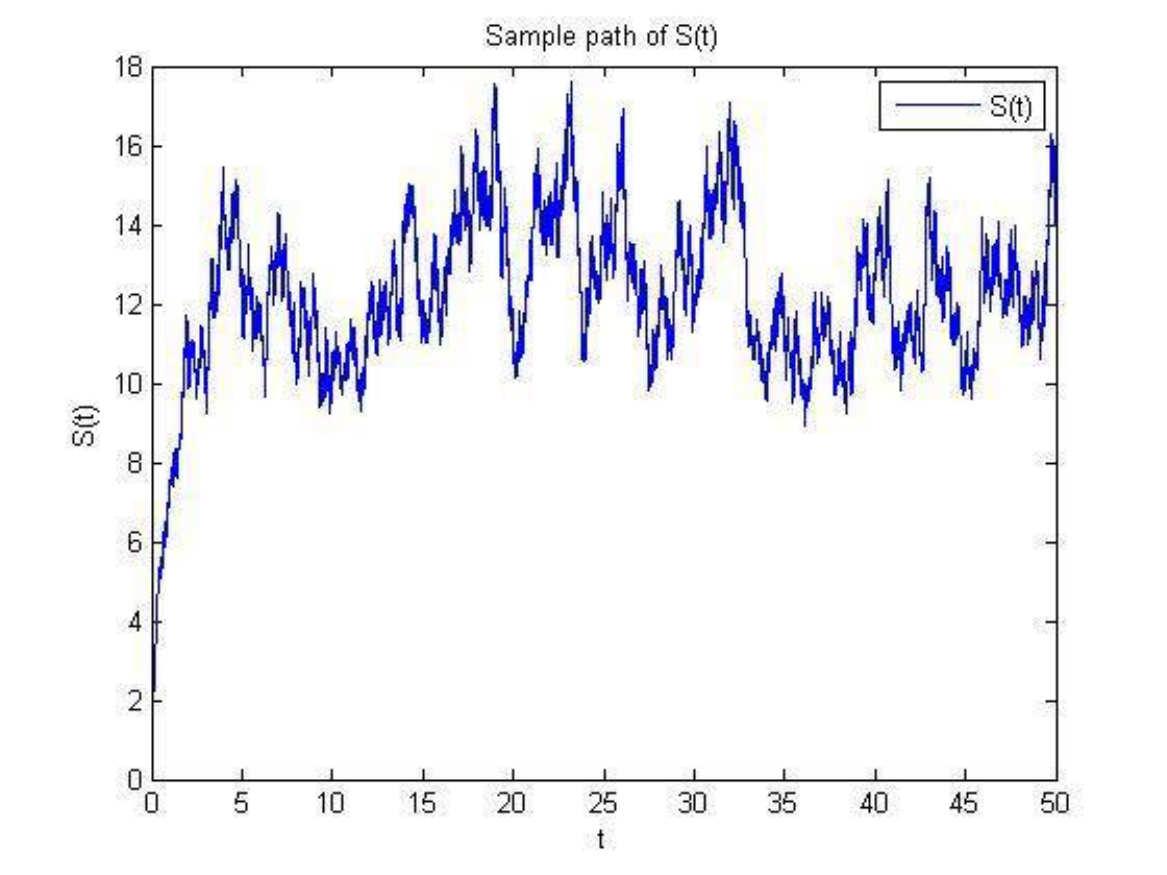}
\includegraphics[width=7.5cm, height=7cm]{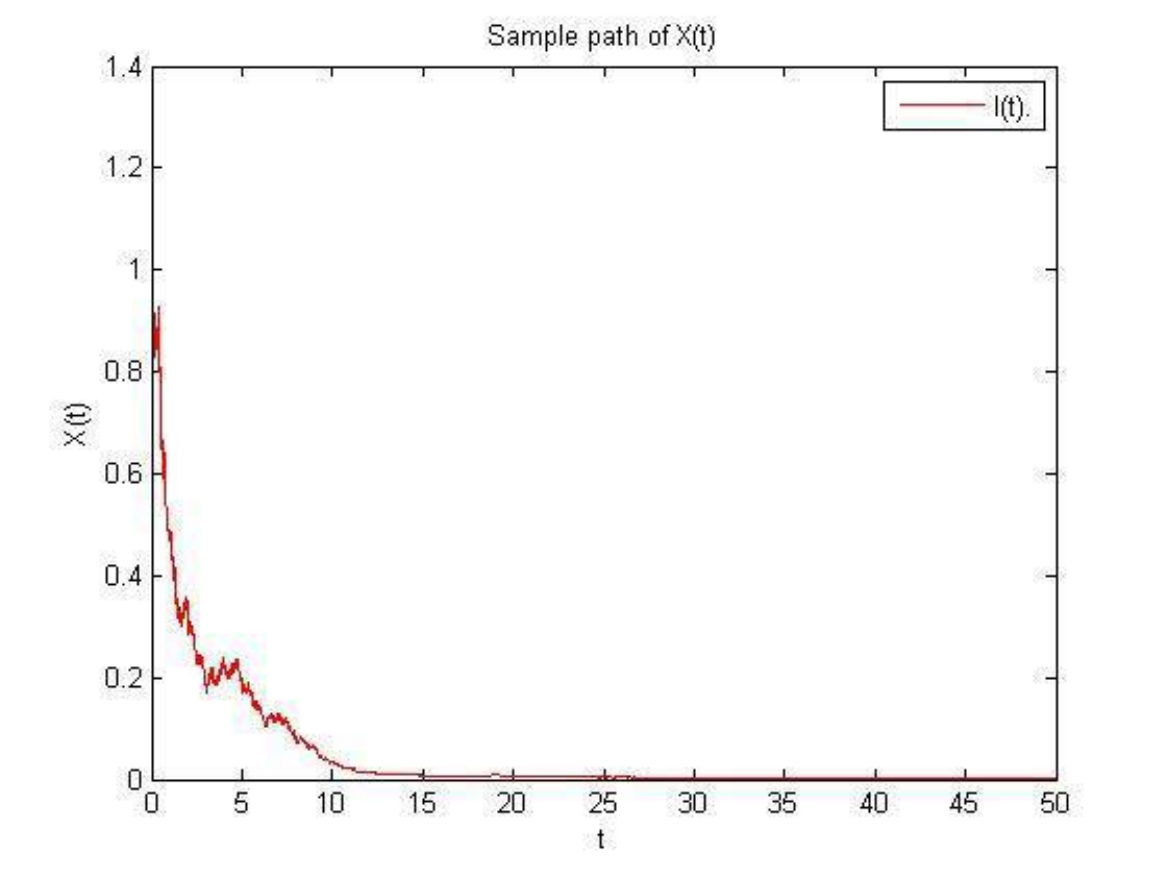}
\caption{Trajectories of $S(t)$ (on the left) and $X(t)$ (on the right) in Example \ref{ex2}.}\label{f2.1}
\end{figure}
}\end{exam}

\begin{exam}\label{ex3}{\rm
Consider \eqref{ww3-111} without switching and parameters
$$
S_0=12,
k_m=8,
\theta=5,
\rho=0,
Y=0.6,
k_d=0.06,
K_S=60,
\sigma_1=0.2, \sigma_2=0.2.
$$
We have $\lambda\approx0.5$.
Sample paths are given in Figure \ref{f3.1}, and
 the density of the empirical measure, which approximate the invariant density,
is shown in Figure \ref{f3.2}.
   \begin{figure}[htp]
\centering
\includegraphics[width=7.5cm, height=7cm]{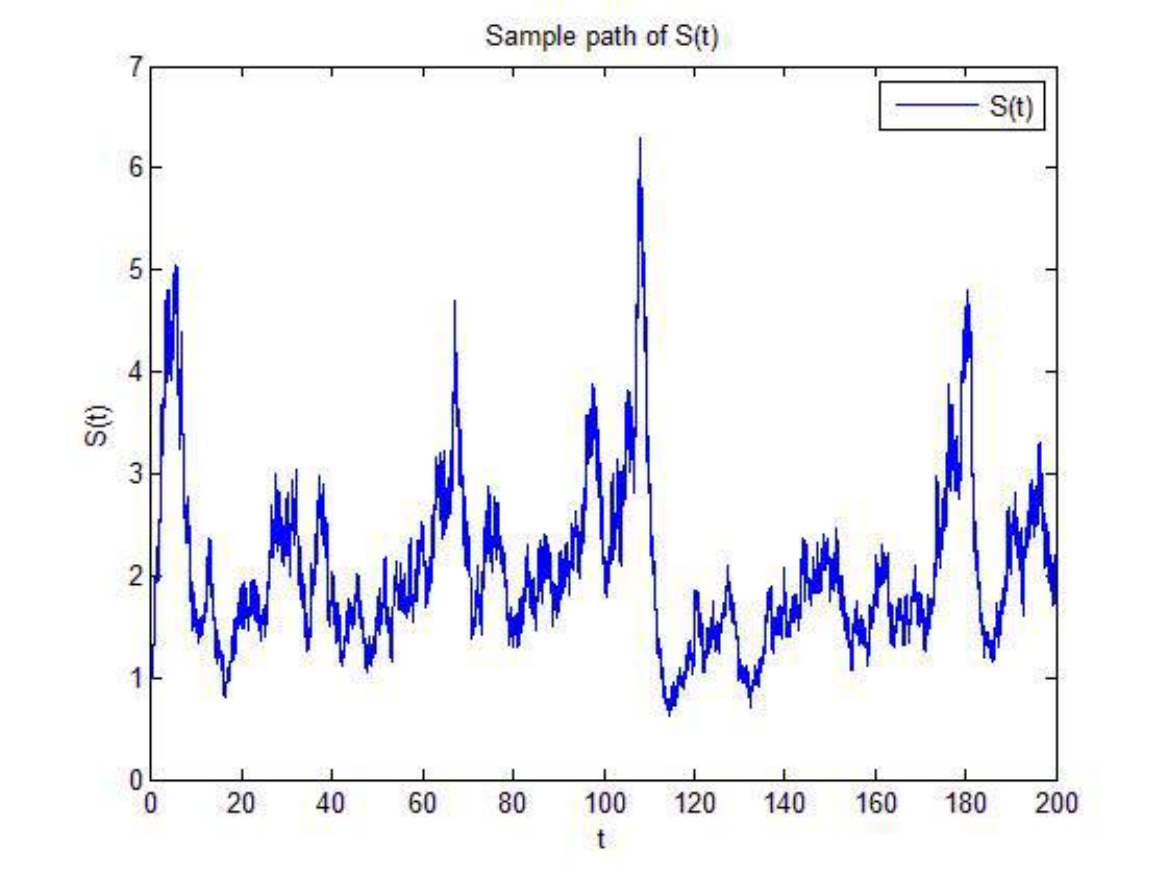}
\includegraphics[width=7.5cm, height=7cm]{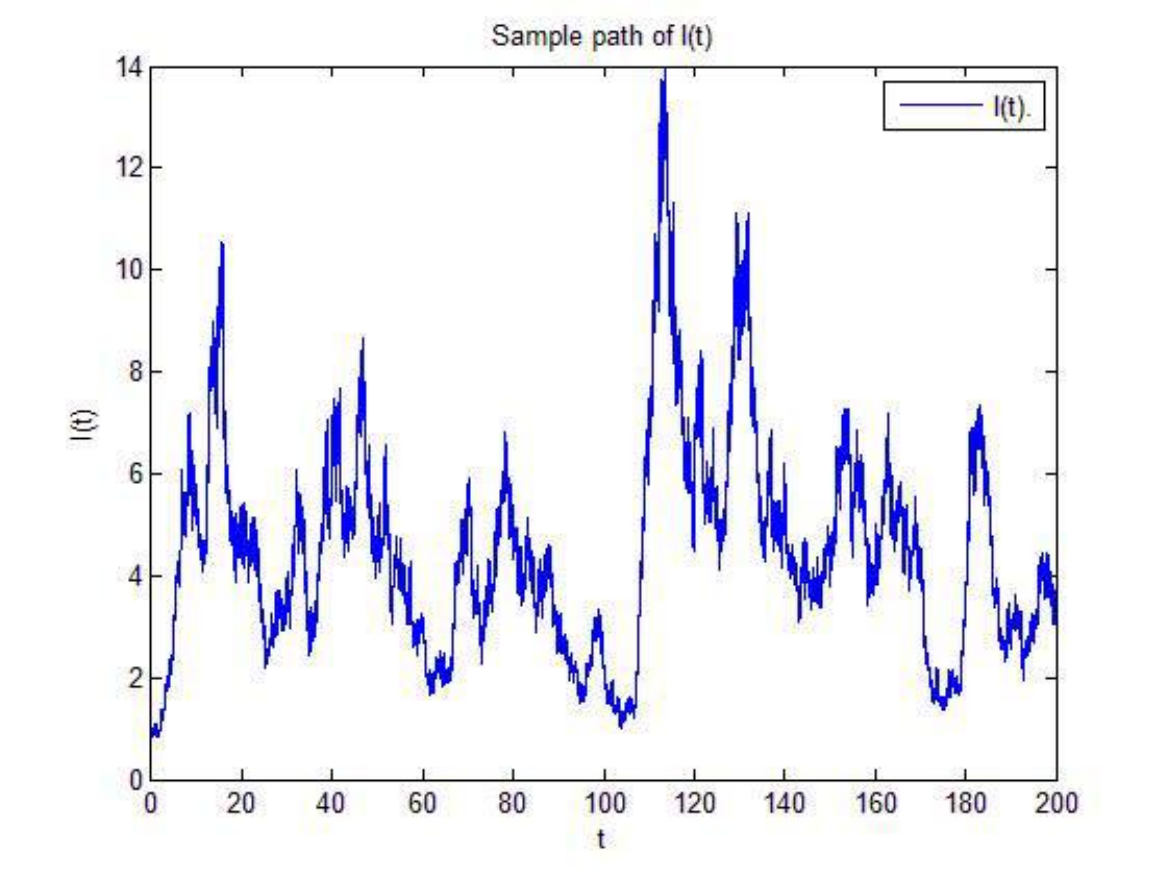}
\caption{Sample paths of $S(t)$ ( on the left) and $X(t)$ (on the right) in Example \ref{ex3}.}\label{f3.1}
\end{figure}
   \begin{figure}[htp]
\centering
\includegraphics[width=7.5cm, height=7cm]{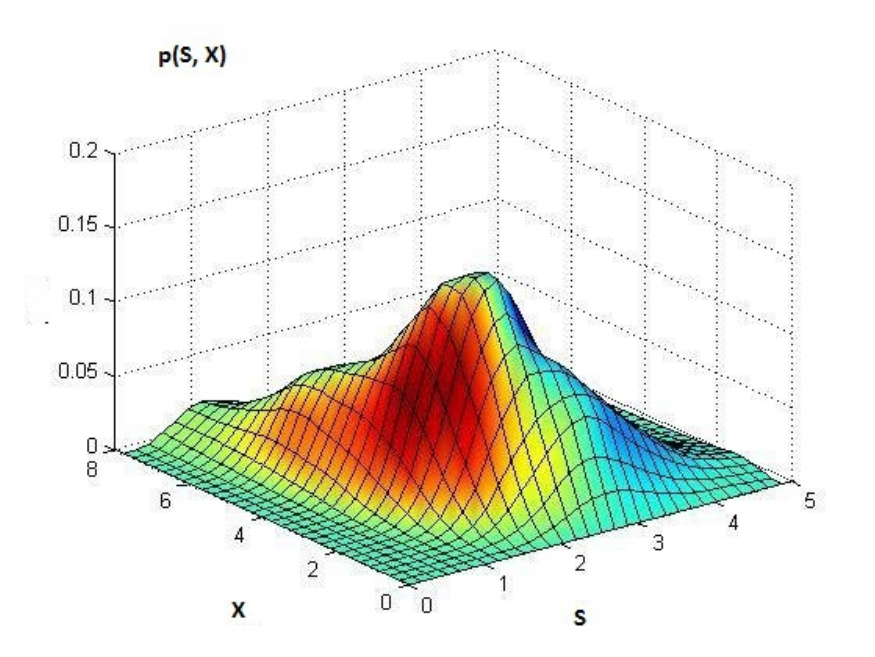}
\includegraphics[width=7.5cm, height=7cm]{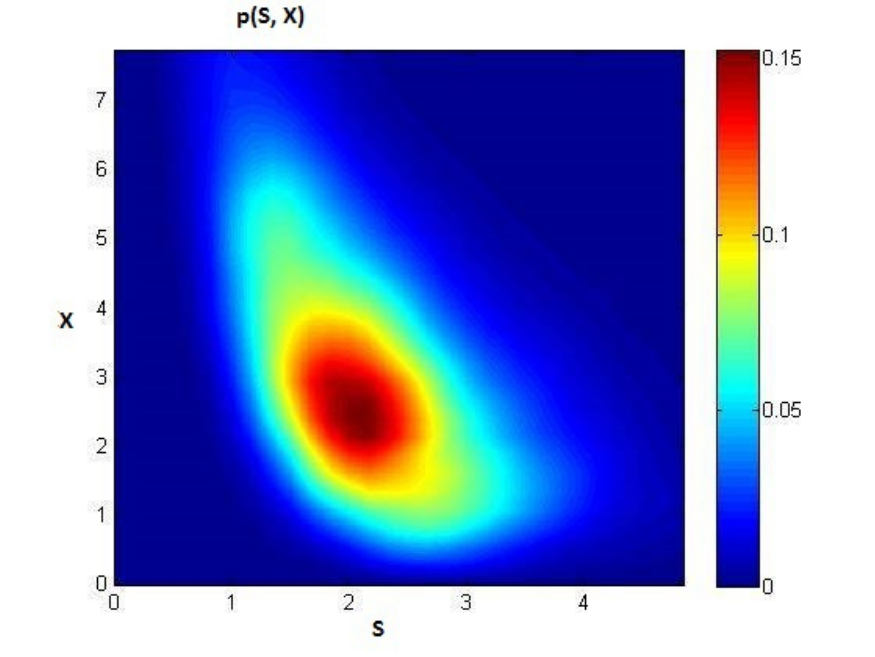}
\caption{The left figure is the 3D graph of the density of the invariant probability measure in Example \ref{ex3}.
The right one depicts the density using scaled colors.}\label{f3.2}
\end{figure}
}
\end{exam}

\begin{exam}\label{ex4}
{\rm The limit $ES^*:=\lim_{t\to\infty}\E_{s,x, i} S(t)$
is regarded as the expected effluent concentration.
We are interested in investigating
 the limit $ES^*$ and $\lambda$
as
 functions of the hydraulic residence time $\theta$.
It can be seen that the expected effluent concentration
is decreasing in $\theta$.
By Theorem \ref{thm2.2} and
\eqref{2-e7},
 we have for $\lambda>0$ that
$$
ES^*:=\lim_{t\to\infty}\E_{s,x,i} S(t)=
\begin{cases}
S_0 &\text{ if } \lambda<0\\
\disp\sum_{j\in\M}\int_{\R^{2,\circ}_+}s\mu^*(ds,dx, j)<S_0&\text{ if }\lambda>0.
\end{cases}
$$
Note that
when $\lambda<0$,  the expected effluent concentration levels off at $S_0$
and then becomes smaller than $S_0$ after $\theta_0$: the value of $\theta$ at which $\lambda=0$.
The numerical approximation (see Figures \ref{f4.1} and \ref{f4.2}) for the expected effluent concentration justifies the claim.
Some fluctuations are due to the errors of
approximation of the random processes.
The behavior of $ES^*$ as a function of $\theta$
is very similar to the deterministic counterpart in \cite{NS}.

When one designs the treatment, a crucial design parameter is the so-called {\it wash-out time}.
If the residence time $\theta$ is less than a critical value, denoted by  $\theta_0$, then the sewage flow is
too fast  for bacteria to grow, existing cells are flushed out faster than they
can multiply.
As a result, the bacteria become extinct.
Figures \ref{f4.1} and \ref{f4.2} show that
$\lambda$ is an increasing function of $\theta$.
By our theoretical results,
to find the wash-out time $\theta_0$,
we need to solve the equation
$\lambda(\theta)=0$.
For the system without switching \eqref{ww2},
the value $\lambda$ can be obtained in a closed form  by solving the Fokker-Planck equation. Then
we can  solve the equation
$\lambda(\theta)=0$ by a standard numerical scheme.
In Figure \ref{f4.2}, we can see that $\theta_0\approx1.4 $.
When random switching are involved, the value of $\lambda$ in \eqref{lambda1}
cannot be solved in a closed form.
However, because of the exponential convergence rate,
one can also perform a numerical approximation to
find $\theta_0$.
In Figure \ref{f4.1}, $\theta_0\approx 0.8$.
      \begin{figure}[htp]
\centering
\includegraphics[width=7.5cm, height=7cm]{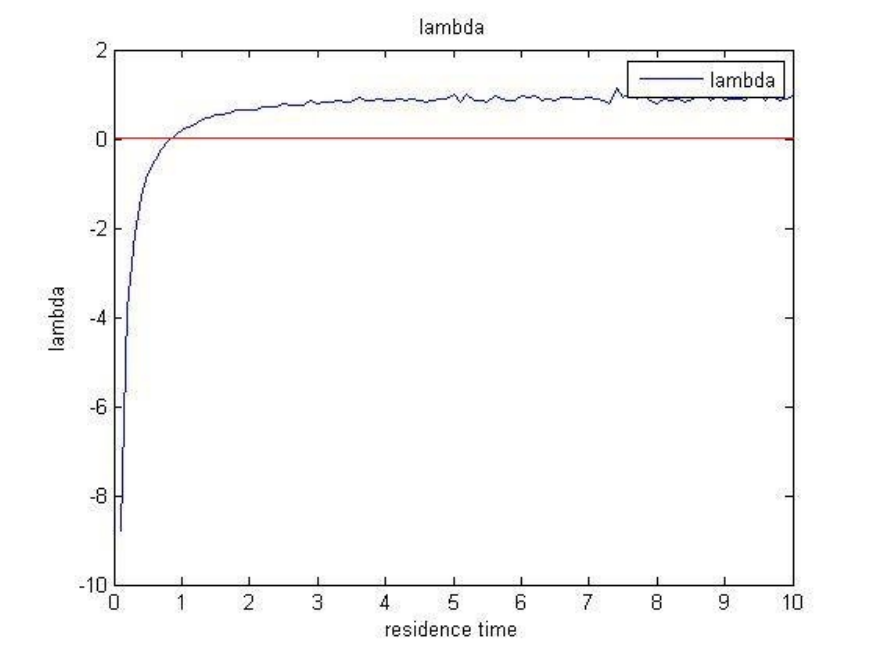}
\includegraphics[width=7.5cm, height=7cm]{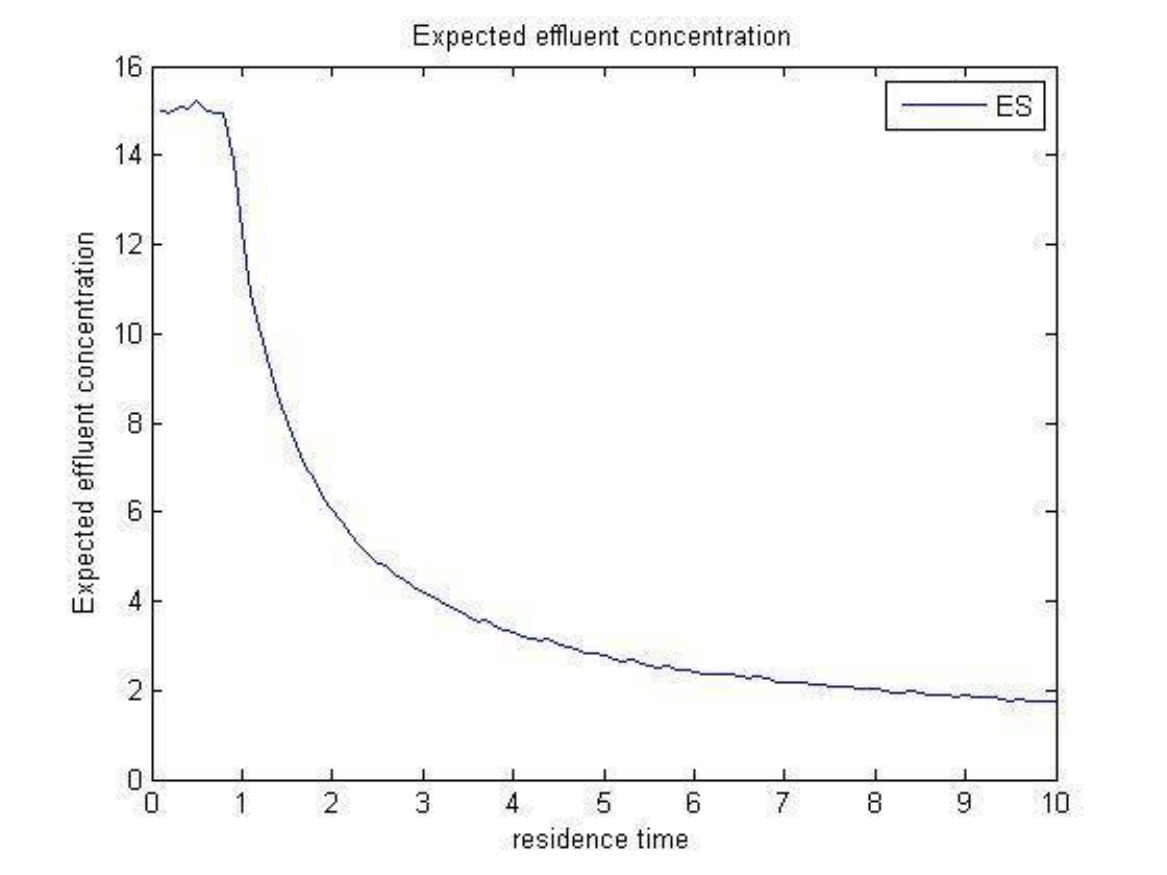}
\caption{$\lambda$ and the
expected effluent concentration as a function of $\theta$ with other parameters as in Example \ref{ex1}.}\label{f4.1}
\end{figure}
   \begin{figure}[htp]
\centering
\includegraphics[width=7.5cm, height=7cm]{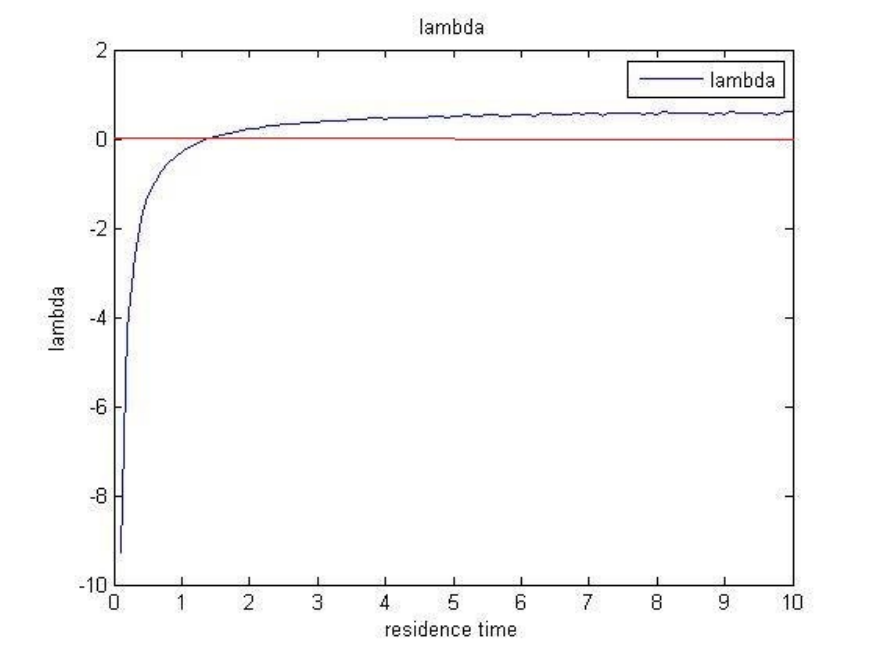}
\includegraphics[width=7.5cm, height=7cm]{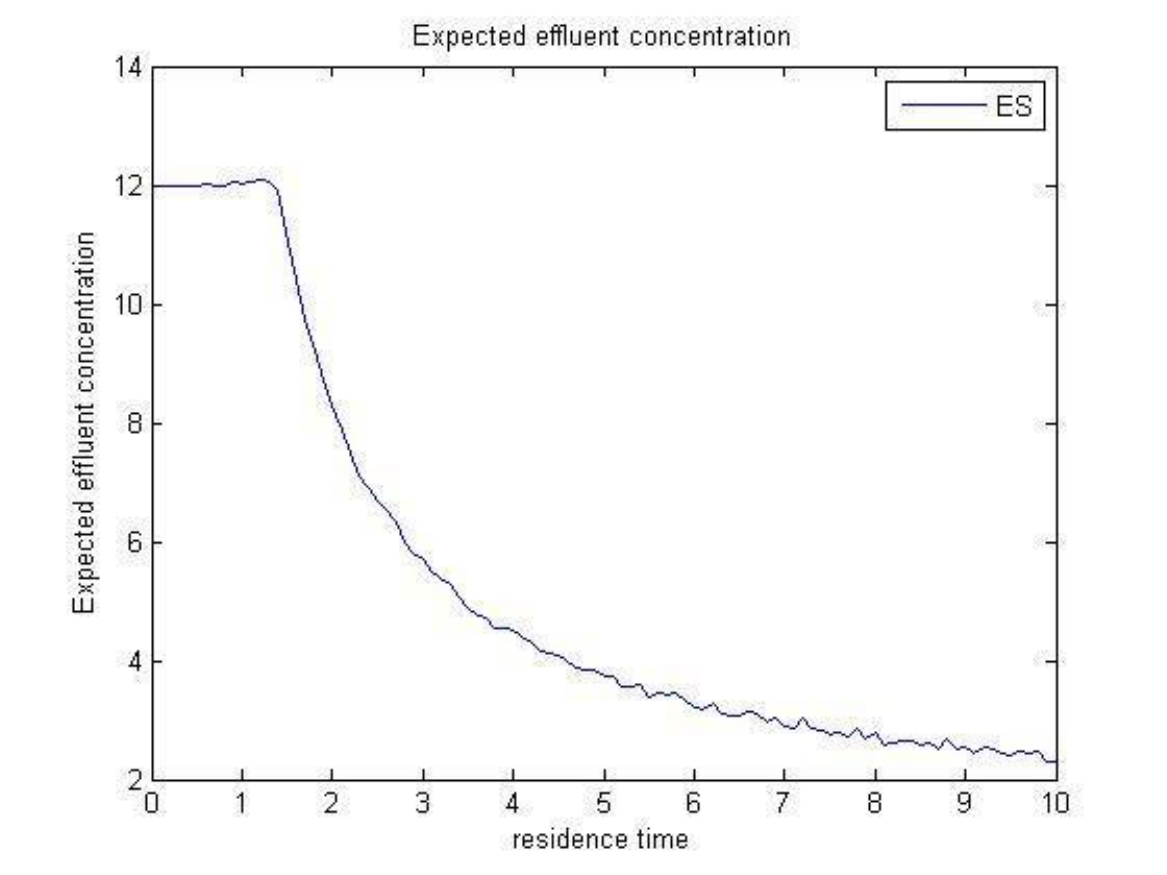}
\caption{$\lambda$ and the
expected effluent concentration as a function of $\theta$ with other parameters as in Example \ref{ex3}.}\label{f4.2}
\end{figure}
}\end{exam}

\section{Controlled Stochastic Chemostat Models}\label{sec:control}
To better help us reaching our goal of dynamically regulating and optimizing the performance of
chemostat models, we introduce a control process in this section and
consider controlled stochastic chemostat models.
Note that for notational simplicity,
we consider only the controlled dynamic systems without switching in this section. Switching can be added, but the notation would be much more complex. One needs to deal with a system of dynamic programming equations  (partial differential equations) in lieu of a single  equation.
It seems to be more instructive to treat relatively simpler models to present the main ideas and leaving out the notational details.

Suppose that we can add a certain amount of bacteria $u(t)$ (the control)
to the system at any time $t$, which
enables us to better adjust the system performance so as to minimize the amount of substrate over a long-time horizon. In this section, we focus on the case $f_0(s,i)=s$.
Then we have a controlled differential equation as follow
\begin{equation}\label{ww-c1}
\begin{cases}
\begin{aligned}
dS(t)=&\left(\dfrac{S_0-S(t)}\theta-X(t)f_1(S(t), X(t))\right)dt+\sigma_1S(t)dW_1(t),\\
dX(t)=&\left(u(t)+X(t)\left(f_2(S(t), X(t))-k_d-\dfrac{1+R}\theta\right)\right)dt+\sigma_2X(t)dW_2(t).
\end{aligned}
\end{cases}
\end{equation}
We assume the control $u(t)$ taking value in a compact interval $[0,M]$ for some $M>0$.
Our objective is to minimize
$$
\limsup_{T\to\infty} \dfrac1T\E_{s,x}^u \int_0^TS(t)dt,
$$
over the class of admissible controls $u(t)$, where $u(t)$ is $\F_t$-adapted. That is, we aim to minimize  the amount of substrate
over the infinite horizon. The cost criterion is in the sense of
an average cost per unit time (or long-run average cost).

Before proceeding further, let us describe how we plan to carry out the analysis. In getting the desired optimal control, we need to obtain the Hamilton-Jacobi-Bellman (HJB) equation.
We use a ``vanishing discount
argument'',
 which
 utilizes some ideas from the book \cite{ABG}. Nevertheless,
in the book, the HJB equation is derived under either ``near-monotone" or ``stable" conditions. Unfortunately,
our model satisfies neither of these conditions.
As a result, some new approaches
are needed to obtain the HJB equation.
The intuitive idea is to look at two different domains.
 In each of  the domains, one of the ``near-monotone" condition
or the ``stable" condition is satisfied. Nevertheless, the complication is that we also need to
 analyze the dynamics of the system, to investigate how the solution moves from one domain to the other, and to
examine how the movement affects the objective function.

To proceed, we provide a road map of our approach.
First, we recall some notation. Then the analysis is carried out using relaxed control setup and appropriate occupation measures.  Theorem \ref{HJB} presents the main result.
To prove it, we need a number of technical results to take care of the situation as mentioned in the last paragraph regarding the two regions satisfying the ``near-monotone" condition or ``stable" conditions separately and the dynamic movements between these regions.
These technical details are presented in a number of lemmas.

To continue,
we recall some concepts and notation introduced in \cite{ABG, KR}.
Let $M(\infty)$
denote the family of measures $\{m(\cdot)\}$ on the Borel subsets of
$[0,\infty)\times [0,M]$
satisfying $m([0,t]\times [0,M])=t$
for all $t\geq0$.
By
the weak convergence $m_n(\cdot)\rightarrow m(\cdot)$
in $M(\infty)$,
we mean $\lim_{n\to\infty}\int f(s,\alpha)m_n(ds\times d\alpha)
=\int f(s,\alpha)m(ds\times d\alpha)
$
for any continuous function $f(\cdot):[0,\infty)\times [0,M]\mapsto\R$ with compact support.
A random measure $m(\cdot)$ with values in $M(\infty)$ is said to be an admissible relaxed control for \eqref{ww-c1}
if $\int_0^M\int_0^tf(s,u)m(ds\times du)$
is independent of $\{W_i(t+s)-W_i(t), s>0, i=1,2\}$ for each bounded and continuous function $f(\cdot)$.
Under a relaxed control $m(\cdot)$, the controlled diffusion \eqref{ww-c1} becomes
\begin{equation}\label{ww-c2}
\begin{cases}
\begin{aligned}
dS(t)=&\left(\dfrac{S_0-S(t)}\theta-X(t)f_1\big(S(t), X(t)\big)\right)dt+\sigma_1S(t)dW_1(t),\\
dX(t)=&\left(\bar m_t+X(t)\left(f_2\big(S(t), X(t)\big)-\wdt k_d\right)\right)dt+\sigma_2X(t)dW_2(t),
\end{aligned}
\end{cases}
\end{equation}
where
$\bar m_t=\int_0^Mu m_t(du)$
and the ``derivative" $m_t$ is defined as the measure-valued function of $(\omega,t)$  such that for any smooth and bounded function $f$, we have
$\iint f(s,u)m(ds\times du)=\int ds \int f(s, u)m_s(du)$.
The operator associated with the controlled diffusion process \eqref{ww-c2}, in which $t$-dependence is hidden,  is given by
$$
\begin{aligned}
\op^m\phi(s,x)=&\dfrac{\partial\phi(s,x)}{\partial s}\left[\frac{S_0-s}\theta-xf_1(s,x)\right]+\dfrac{\partial\phi(s,x)}{\partial x} \left[\bar m_t+xf_2(s,x)-x\wdt k_d\right]\\
&+\dfrac12\left(\sigma_1^2\dfrac{\partial^2\phi(s,x)}{\partial s^2}s^2+\sigma_2^2\dfrac{\partial^2\phi(s,x)}{\partial x^2}x^2\right).
\end{aligned}
$$

\begin{defn}{\rm We have the following definitions and
notations.
\begin{itemize}
\item  Let $\mathcal P(M(\infty))$ be the space of probability measures on $M(\infty)$. A relaxed control $m(\cdot)$ for \eqref{ww-c2} is said to be Markov
if there exists a measurable function $v:\R^2_+\mapsto\mathcal P(M(\infty))$
such that $m_t=v(S(t), X(t)), t\geq0.$
Under a relaxed Markov control $m_t=v(S(t),X(t))$, the solution process $(S(t), X(t))$ to \eqref{ww-c2} is a Markov process with generator
$$
\begin{aligned}
\op^v\phi(s,x)=&\dfrac{\partial\phi(s,x)}{\partial s}\left[\frac{S_0-s}\theta-xf_1(s,x)\right]+\dfrac{\partial\phi(s,x)}{\partial x} \left[v(s,x)+xf_2(s,x)-x\wdt k_d\right]\\
&+\dfrac12\left(\sigma_1^2\dfrac{\partial^2\phi(s,x)}{\partial s^2}s^2+\sigma_2^2\dfrac{\partial^2\phi(s,x)}{\partial x^2}x^2\right).
\end{aligned}
$$
\item A Markov control $v$ is a relaxed control satisfying that
$v(z)$ is a Dirac measure on $[0,M]$ for each $z\in\R^2_+$.
\item Denote the set of Markov controls
and relaxed Markov controls by $\Pi_{M}$ an $\Pi_{RM}$, respectively.
With a relaxed Markov control, $(S(t), X(t))$
is a Markov process that has the strong Feller property in $\R^{2,\circ}_+$;
see \cite[Theorem 2.2.12]{ABG}.
\item Since the diffusion is nondegenerate in $\R^{2,\circ}_+$,
if the process $(S(t), X(t))$ has an invariant probability measure in $\R^{2,\circ}_+$,
the invariant measure is unique, denoted by $\eta_v$.
In this case, the control $v$ is said to be stable.
Denote by $\Pi_{SRM}$ the set of stable relaxed Markov controls.
\item Let $\mathcal P(\mathcal X)$ be the space of probability measures on a metric space $\mathcal X$. For any stable relaxed Markov control $v$,
define
$$\pi_v(dz\times du)=[v(z)(du)]\times \eta_v(du)\in \mathcal P(\R^{2,\circ}_+\times [0,M]),$$
and
$$\mathcal{G}=\{\pi_v: v \text{ is a stable relaxed Markov control }\}\subset\mathcal P(\R^{2,\circ}_+\times [0,M]).$$
\end{itemize}
}
\end{defn}	

We need the following lemma whose proof is analogous to \cite[Lemma 2.3]{DNDY}.

\begin{lm}\label{lm3.1}
There exist a  sufficiently small $p>0$ and positive constants $K_1$, $K_2$, and $K_3$ such that
\begin{equation}\label{lm3.1-e0}
\op^m (2\kappa_0 s+x)\leq K_1-K_2 (2\kappa_0s+x),
\end{equation}
\begin{equation}\label{lm3.1-e5}
\op^m (2\kappa_0 s+x)^{1+p}\leq K_1-K_2 (2\kappa_0s+x)^{1+p}-K_2 (2\kappa_0s+x)^p xf_1(s,x),
\end{equation}
and
\begin{equation}\label{lm3.1-e1}
\op^m (2\kappa_0 s+x)^4\leq K_3(1+2\kappa_0s+x)^4,
\end{equation}
for any admissible relax control $m(\cdot)$.
Consequently, it holds for any admissible relaxed control $m(\cdot)$ that
\begin{equation}\label{lm3.1-e2}
   \E_{s,x}^m\big(2\kappa_0S(t)+ X(t)\big)^{1+p} \leq (2\kappa_0s+x)^{1+p}e^{-K_2t}
  +\frac{K_2}{K_1}\text{ for } s>0, x\geq 0,
\end{equation}
\begin{equation}\label{lm3.1-e3}
\E^m_{s,x} \int_0^T\big(2\kappa_0S(t)+X(t)\big)^p X(t)f_1(S(t),X(t))dt\leq (2\kappa_0s+x)^{1+p} +K_1 T,
\end{equation}
and
\begin{equation}\label{lm3.1-e6}
   \E_{s,x}^m\big(2\kappa_0S(t)+ X(t)\big)^{4} \leq (1+2\kappa_0s+x)^{4}e^{K_3t}
  \text{ for } s>0, x\geq 0, t\geq0.
\end{equation}
As a result, we have
\begin{equation}\label{lm3.1-e4}\barray
\ad
\int_{\R^{2,\circ}_+}(2\kappa_0s+x)^{1+p}\eta_v(ds,dx)\leq \frac{K_2}{K_1} \text{ and } \\
 \ad \int_{\R^{2,\circ}_+}(2\kappa_0s+x)^p xf_1(s,x)\eta_v(ds,dx)\leq \frac{K_2}{K_1}\,\text{ for all } v\in\Pi_{RM}.\earray
\end{equation}

\end{lm}
Define the following sets
\begin{equation}\label{defh0}
\h_0:=\left\{(s,x)\in\R^{2,\circ}_+: 2\kappa_0s+x\leq \frac{K_1+1}{K_2}+2\kappa_0:=K_4 \right\},
\end{equation}
and
\begin{equation*}
\h=\left\{(s,x)\in\R^{2,\circ}_+: 2\kappa_0s+x\leq K_4
\text{ and } s\geq \delta_0\right\}\subset\h_0,
\end{equation*}
where $\delta_0$ is a positive constant
to be determined in the proof of Lemma 4.2.
For a closed set $\mathcal K\in \R^{2,\circ}_+,$
define $\tau_{\mathcal K} =\inf\{t\geq 0: (S(t), X(t))\in \mathcal K\}$.
Since
$\op^u (2\kappa_0s+x)\leq -1$ for $(s,x)\notin\h_0, u\in[0,M]$,
we have
\begin{equation}\label{e3.2}
\E^m_{s,x}\tau_{\h_0}\leq 2\kappa_0s+x \,\text{ for any admissible relaxed control } \, m(t).
\end{equation}
We have the following lemmas whose proofs are given in the appendix.

\begin{lm}\label{lm3.3}
There is a constant $C_1>0$ depending only on $\h$ and $\h_0$ such that
$$\E^v_{s,x} \tau_\h\leq C_1\,\text{ for }\, (s,x)\in\h_0, v\in\Pi_{RM}.$$
Moreover,
\begin{equation}\label{lm3.3-e0}
\E^v_{s,x} \tau_\h\leq 2\kappa_0s+x+C_1, (s,x)\in\R^{2,\circ}_+, v\in\Pi_{RM}.
\end{equation}
With the constant control $v_c\equiv M$, we have
\begin{equation}\label{lm3.3-e0a}
\E^{v_c}_{s,x} \tau_{\h_*}\leq 2\kappa_0s+x+\hat C_1, (s,x)\in\R^{2,\circ}_+,
\end{equation}
where $\h_*$ is a compact subset of $\h_0$, and $\hat C_1$ is a positive constant depending on $\h_0$ and $\h_*$.
\end{lm}

\begin{lm}\label{lm3.6}
For any $L_1>0$ and $\eps>0$, there exists a $\delta>0$ such that
\begin{equation}\label{lm3.6-e0}
\PP_{s,\eps}^v\{X(\tau_\h)\geq\delta\}\geq \dfrac 12,\text{ for any }s\leq L_1, v\in \Pi_{RM}.
\end{equation}
\end{lm}

With these lemmas, let
$$\rho_v:=\int_{\R^{2,\circ}_+\times [0,M]} s\pi_v(ds\times dx\times du) \text{ and } \rho^*=\inf_{v\in\Pi_{SRM}}\left\{\rho_v\right\}.$$
Since \eqref{lm3.3-e0a} implies
the existence of an invariant probability measure for $(S(t), X(t))$ under control $v_c$, we claim that $\Pi_{SRM}\ne\emptyset$.
Moreover, for any admissible relaxed control $m(t)$,
we have that
$$\dfrac{\E_{s,x}^{m}S(T)-S(0)}T= \E_{s,x}^{m}\dfrac1T\int_0^T\left(\dfrac{S_0-S(t)}\theta- X(t)f_1\big(S(t), X(t)\big)\right)dt.$$
In view of \eqref{lm3.1-e2}, we have
\begin{equation}\label{Sinfty}
\lim_{T\to\infty}\E_{s,x}^{m}\dfrac1T\int_0^T
\left(\dfrac{S_0-S(t)}\theta-X(t)f_1\big(S(t), X(t)\big)\right)dt=\liminf_{T\to\infty} \dfrac{\E_{s,x}^{m}S(T)-S(0)}T=0,
\end{equation}
which leads to
\begin{equation}\label{LSeq}
\int_{\R^{2,\circ}_+}\left(\dfrac{S_0-s'}\theta-x'f_1(s',x')\right)\eta_{v}(ds',dx')=0 \,\text{ for any }\,v\in\Pi_{SRM}.
\end{equation}
As a result,
\begin{equation}\label{rhos0}
\begin{aligned}
\rho^*\leq \lim_{T\to\infty}\E_{s,x}^{v_c}\dfrac1T\int_0^TS(t)dt
=&\int_{\R^{2,\circ}_+}s'\eta_{v_c}(ds',dx')\\
=&\theta\left[\dfrac{S_0}\theta - \int_{\R^{2,\circ}_+}x'f_1(s',x')\eta_{v_c}(ds',dx')\right]\\
<&S_0.
\end{aligned}
\end{equation}

To proceed,
we  derive a lemma, which allows us to find an optimal control
in $\Pi_{SRM}$.

\begin{lm}\label{lm3.2}
For any admissible relaxed control $m$, define an empirical
measure
$\zeta_T^m$  as a $\PU$-valued process
satisfying that
$$
\int_{\UU}fd\zeta_T^m=\E^m_{s,x}\frac1T\int_0^T\left(\int_0^Mf(S_t, X_t, u)m_t(du)\right)dt.
$$
Then,
with probability 1, every limit point of $\zeta_T^m$ as $T\to\infty$ can be decomposed as
\begin{equation}\label{lm3.2-e1}
\hat\zeta=\delta\zeta'+(1-\delta)\zeta'',
\end{equation}
where $\zeta'\in\mathcal{G}$
and $\zeta''\in \mathcal{P}((0,\infty)\times\{0\}\times[0,M])$
satisfying
$$\int_{(0,\infty)\times\{0\}\times[0,M]}s\zeta''(ds,dx,du)=S_0.$$
As a result,
for any admissible relaxed control $m$,
\begin{equation}\label{lm3.2-e2}\liminf_{T\to\infty} \dfrac1T\E_{s,x}^m\int_0^TS(t)dt\geq \rho^*.
\end{equation}
\end{lm}

Lemma \ref{lm3.2} enables us to find an optimal control that is a Markov relaxed control.
To find the Markov control, we need to establish the HJB equation associated to the control problem \eqref{ww-c2}.
Let $\mathcal{C}_{1+p}$ be the class of functions $V: \R^{2,\circ}_+\mapsto\R$
such that $V(s,x)\leq c_V(1+s+x)^{1+p'}, (s,x)\in\R^{2,\circ}_+$ for some constants $c_V>0$
and $p'\in(0,p)$.
The rest of this section aims to prove a theorem on the existence and uniqueness of
solutions to the HJB equation.

\begin{thm}\label{HJB}
There is a unique pair
$(V,\rho)$, where $V\in C^2(\R^{2,\circ}_+)\cap\mathcal{C}_{1+p}$
and $\rho\in\R$ satisfying the equation
$$\min_{u\in[0,M]}\left\{\op^u V(s, x)+s\right\}=\rho.$$
Moreover, we have $\rho=\rho^*$
and $v^*\in\Pi_{RM}$ is an optimal control if and only if
it is a measurable selector from the minimizer
$$\min_{u\in[0,M]}\left\{\op^u V(s, x)+s\right\}.$$
In fact, we can choose
$$
v^*(s,x)=\begin{cases}
0&\,\text{ if }\, \dfrac{\partial V(s,x)}{\partial x}\geq 0,\\
M&\,\text{ otherwise. }
\end{cases}
$$
\end{thm}

In \cite{ABG}, the HJB equation is obtained under either of the
 following two assumptions:
(a) the cost function satisfies the so-called {\it near-monotone condition}, or
(b) any relax Markov control is stable and the set $\{\eta_v, v\in\Pi_{RM}\}$ is tight.
However, neither of these conditions are satisfied for our system.
Thus the results in \cite{ABG} are not directly applicable.

To prove the existence and uniqueness of solutions to the HJB equation, we use an idea that might be called a vanishing discount argument. That is,
we examine
the cost and value functions
of the corresponding discounted control problem
and look at the limit when the discount factor tends to $0.$
We also need to estimate the value function in different parts of $\R^{2,\circ}_+$ to obtain desired properties of the value functions, which is key to prove Theorem \ref{HJB}.
That is the most difficult task of this section.

Let $V_\gamma(s,x)$ be the optimal $\gamma$-discounted cost, that is
$$V_\gamma(s,x)=\inf\left\{\E^m_{s,x}\int_0^\infty e^{-\gamma t} S(t): m \text{ runs over the set of relaxed controls }\right\}, \, (s,x)\in\R^{2,\circ}_+.$$
Then it follows from \cite[Theorem 3.5.6 \& Remark 3.5.8]{ABG} that $V_\gamma(s,x)\in C^2(\R^{2,\circ}_+)\cup C_b(\R^{2,\circ}_+)$ satisfies
\begin{equation}\label{e3.5}
\min_{u\in[0,M]}\left\{\op^uV_\gamma(s,x)+s\right\}=\gamma V_\gamma(s,x).
\end{equation}
and the optimal Markov control $v_\gamma$ is a selector
of $\min_{u\in[0,M]}\left\{\op^uV_\gamma(s,x)+s\right\}$.
The following lemma is
from \cite[Lemma 3.7.8]{ABG}.

\begin{lm}\label{lm3.4}
Fix $(s_*,x_*)\in\R^{2,\circ}_+$.
For any sequence $\gamma_n\downarrow 0$,
there exists a subsequence, which is still denoted by $\{\gamma_n\}$, and a function $V\in C(\R^{2,\circ}_+)$ and a constant $\rho$ such that
as $n\to\infty$, we have
\begin{equation}\label{lm3.4-e0}
\gamma_n V_{\gamma_n}(s_*, x_*)\to \rho\,\text{ and }\, \bar V_{\gamma_n}(s,x):=V_{\gamma_n}(s,x)-V_{\gamma_n}(s_*,x_*)\to V(s,x)
\end{equation}
uniformly on each compact subset of $\R^{2,\circ}_+$.
Moreover, we have
$$\min_{u\in[0,M]}\left\{\op^u V(s, x)+s\right\}=\rho\leq\rho^*, \, (s,x)\in\R^{2,\circ}_+.$$
\end{lm}

We aim to show that a limit function $V$ in Lemma \ref{lm3.4}
belongs to the family $\mathcal{C}_{1+p}$.

\begin{prop}\label{lm3.5}
Let $V$ be any limit in \eqref{lm3.4-e0}. Then
$$\sup_{(s,x)\in\h}|V(s,x)|<\infty.$$
\end{prop}

\begin{proof}
Denoted by $v_\gamma$ the optimal Markov control of the $\gamma$-discounted control problem.
Let $\eps_*=\frac{S_0-\rho^*}5$,
and $\bar\rho=\rho^*+\eps_*=S_0-4\eps^*$ and $C_1$ be as in Lemma \ref{lm3.3}.
By \eqref{lm3.1-e2},
there exists a $C_2>0$ depending only on $\h_0$ such that
\begin{equation}\label{bd-S}
\E^{v_\gamma}_{s,x}(2\kappa_0S(t)+X(t))\leq C_2, \text{ for any } (s,x)\in\h_0, v\in\Pi_{RM}, t\geq 0.
\end{equation}
By It\^o's formula,
we have that
$$\E^{v_\gamma}_{s,x}\left[S(t+h)-S(t)\right]=\int_t^{t+h}\left[\dfrac{S_0-\E^{v_\gamma}_{s,x}S(y)}\theta-\E^{v_\gamma}_{s,x}X(y)f_1\big(S(y), X(y)\big)\right]dy.
$$
Letting $h\to0$, we obtain
\begin{equation}\label{lm3.5-e2}
\begin{aligned}
\dfrac{d}{dt}\E^{v_\gamma}_{s,x}S(t)=&\dfrac{S_0-\E^{v_\gamma}_{s,x}S(t)}\theta-\E^{v_\gamma}_{s,x}X(t)f_1\big(S(t), X(t)\big).
\end{aligned}
\end{equation}
Consider the differential equation
\begin{equation}\label{e3.10}
dz(t)=\frac{S_0-\eps_*-z}\theta dt.
\end{equation}
It is easy to show that $z(t)$ converges to $S_0-\eps_*$
uniformly for each initial value belonging to any bounded set.
Thus, there is a $T^*>\frac{\bar\rho(C_1+C_2)+1}{\eps_*}$ such that
\begin{equation}\label{lm3.5-e3}
\int_0^{T^*} z(t)dt\geq (S_0-2\eps_*)T^*, \text{ for any } z(0)\geq 0.
\end{equation}
In view of \eqref{lm3.4-e0},
we can choose $\gamma_*$ sufficiently small such that
\begin{equation}\label{lm3.5-e6}
(1-e^{-\gamma_* T^*})C_2\leq \eps_*,
\end{equation}
and
\begin{equation}\label{lm3.5-e5}
\gamma V_{\gamma}(s_*,x_*)\leq\bar\rho, \,\text{ for }\, \gamma<\gamma_*.
\end{equation}
We divide $\h$ into two subsets (one of which can be empty):
$$
\h^\gamma_1=\left\{(s,x)\in\h: \E^{v_\gamma}_{s,x} X(t)f_1\big(S(t), X(t)\big) \leq\frac{\eps_*}{\theta}, \text{ for all } t\in[0,T^*]\right\}\, \text{ and }\, \h^\gamma_2=\h\setminus \h^\gamma_1.
$$

{\bf Step 1:}
Consider $(s,x)\in\h^\gamma_1$.
It follows from \eqref{lm3.5-e2}, \eqref{lm3.5-e3}, and a comparison argument of
differential inequalities that
\begin{equation}\label{lm3.5-e4}
\E_{s,x}^{v_\gamma}\int_0^{T^*} S(t)dt\geq (S_0-2\eps_*)T^*, (s,x)\in\h_1^\gamma.
\end{equation}
Then, we derive from \eqref{bd-S} and \eqref{lm3.5-e6} that
\begin{equation}\label{lm3.5-e1111}
\E_{s,x}^{v_\gamma}\int_0^{T^*} (1-e^{-\gamma t})S(t)dt\leq \eps_*T^*,\, (s,x)\in \h^\gamma_1, \gamma<\gamma_*.
\end{equation}
As a result of \eqref{lm3.5-e4} and \eqref{lm3.5-e1111}
\begin{equation}\label{lm3.5-e7}
\E_{s,x}^{v_\gamma}\int_0^{T^*} e^{-\gamma t}S(t)dt\geq (S_0-3\eps_*)T^*,\,\text{ for }\, (s,x)\in \h^\gamma_1, \gamma<\gamma_*.
\end{equation}
In view of \eqref{lm3.5-e5},
we have that
$$0\leq\inf_{(s,x)\in\h} V_\gamma(s,x)\leq\dfrac{\bar\rho}\gamma<\infty\,\text{ for }\, \gamma\leq \gamma_*,$$
which combined
with the strong Markov property of $(S(t), X(t))$ under a Markov control,
implies
\begin{equation}\label{lm3.5-e8}
\begin{aligned}
V_\gamma(s,x)=&\E_{s,x}^{v_\gamma}\int_0^{\tau^{T^*}_\h} e^{-\gamma t}S(t)dt + \E_{s,x}^{v_\gamma} \left[e^{-\gamma \tau^{T^*}_\h} V_\gamma\left(S(\tau^{T^*}_\h), X(\tau^{T^*}_\h)\right)\right]\\
\geq&
\E_{s,x}^{v_\gamma}\int_0^{\tau^{T^*}_\h} e^{-\gamma t}S(t)dt + \left(\inf_{(s,x)\in\h} V_\gamma(s,x)\right)\E_{s,x}^{v_\gamma} e^{-\gamma \tau^{T^*}_\h} \\
\geq&\E_{s,x}^{v_\gamma}\int_0^{\tau^{T^*}_\h} e^{-\gamma t}S(t)dt + \inf_{(s,x)\in\h} V_\gamma(s,x) -\dfrac{\bar\rho}\gamma\E_{s,x}^{v_\gamma} (1-e^{-\gamma\tau^{T^*}_\h})\,\text{ for }\, (s,x)\in\R^{2,\circ}_+,\\
\end{aligned}
\end{equation}
where
$\tau^{T^*}_\h=\inf\{t\geq T^*: (S(t),X(t))\in\h\}$.
By the Markov property of $(S(t), X(t))$, \eqref{bd-S}, and \eqref{lm3.3-e0}, we have
$$
\begin{aligned}
\dfrac1{\gamma}\E^{v_\gamma}_{s,x} (1-e^{-\gamma\tau^{T^*}_\h})\leq&\E^{v_\gamma}_{s,x}\tau^{T^*}_\h
\leq T^*+\E^{v_\gamma}_{s,x}(\tau^{T^*}_\h-T^*)\\
=&T^*+\E^{v_\gamma}_{s,x}\E^{v_\gamma}_{S(T^*), X(T^*)}\tau_\h\\
\leq& T^*+C_1+C_2.
\end{aligned}
$$
This together with \eqref{lm3.5-e7} and \eqref{lm3.5-e8} implies
$$
\begin{aligned}
V_\gamma(s,x)\geq& \inf_{(s,x)\in\h} V_\gamma(s,x)+ (S_0-3\eps_*)T^*-\bar\rho(T^*+C_1+C_2)\\
\geq&\inf_{(s,x)\in\h} V_\gamma(s,x)+ \eps_*T^*-\bar\rho(C_1+C_2)\,\text{ (since } \bar\rho=S_0-4\eps_*)\\
> &\inf_{(s,x)\in\h} V_\gamma(s,x)+1\,\text{ (since } T^*>\dfrac{\bar\rho(C_1+C_2)+1}{\eps_*}),
\end{aligned}
$$
which leads to
$$
\bar V_\gamma(s,x)\geq \inf_{(s,x)\in\h} \bar V_\gamma(s,x) +1, (s,x)\in\h^\gamma_1, \gamma\leq\gamma_*.
$$
As a result,
\begin{equation}\label{lm3.5-e9}
\inf_{(s,x)\in\h} \bar V_\gamma(s,x)= \inf_{(s,x)\in\h^\gamma_2} \bar V_\gamma(s,x) \,\text{ for }\gamma\leq\gamma_*.
\end{equation}

{\bf Step 2:} Consider $(s,x)\in\h_2^\gamma$.
In view of \eqref{lm3.1-e6}, there is a $C_3>0$ depending only on $\h$ and $T^*$ such that
\begin{equation}\label{lm3.5-e10}
\E^{v_\gamma}_{s,x}(S^4(t)+X^{4}(t))\leq C_3, \text{ for any } t\in[0, T^*], (s,x)\in\h, v\in\Pi_{RM}.
\end{equation}
For $(s,x)\in\h_2^\gamma$, there is a $t_{s,x}\in [0,T^*]$ such that
$\E^{v_\gamma}_{s,x}X(t_{s,x})f_1(S(t_{s,x}), X(t_{s,x}))\geq\dfrac{\eps_*}{\theta}$.
Since $f_1(s,x)\leq\kappa_0 (s+1)$,
$$\left[\E^{v_\gamma}_{s,x}X(t_{s,x})f_1(S(t_{s,x}), X(t_{s,x}))\right]^2\leq \kappa_0^2\E^{v_\gamma}_{s,x}X^2(t_{s,x})\E^{v_\gamma}_{s,x} (S(t_{s,x})+1)^2,
$$
which together with \eqref{lm3.5-e10} implies
$$\E^{v_\gamma}_{s,x}X^2(t_{s,x})\geq \hat\eps_{*},$$
for some $\hat\eps_{*}\in(0,\frac12)$ depending on $\eps_*$ and $C_3$.
Then
$$\E^{v_\gamma}_{s,x}\left[\1_{\{X^2(t_{s,x})<\frac{\hat\eps_*}2\}} X^2(t_{s,x})\right]+\E^{v_\gamma}_{s,x}\left[\1_{\{X^2(t_{s,x})\geq \frac{\hat\eps_*}2\}}X^2(t_{s,x})\right]\geq \hat\eps_*,$$
which leads to
$$\E^{v_\gamma}_{s,x}\left[\1_{\{X^2(t_{s,x})\geq \frac{\hat\eps_*}2\}}X^2(t_{s,x})\right]\geq \dfrac{\hat\eps_*}{2}.$$
By H\"older's inequality and \eqref{lm3.5-e10}, we have
$$
\begin{aligned}
\left(\E^{v_\gamma}_{s,x}\left[\1_{\{X^2(t_{s,x})\geq \frac{\hat\eps_*}2\}}X^2(t_{s,x})\right]\right)^2
\leq& \PP^{v_\gamma}_{s,x}\left\{X^2(t_{s,x})\geq \frac{\hat\eps_*}2\right\}\left[\E^{v_\gamma}_{s,x}X^{4}(t_{s,x})\right]\\
\leq&C_3\PP^{v_\gamma}_{s,x}\left\{X^{2}(t_{s,x})\geq \dfrac{\hat\eps_*}2\right\}.
\end{aligned}
$$
Thus,
$$\PP^{v_\gamma}_{s,x}\left\{X^{2}(t_{s,x})\geq \frac{\hat\eps_*}{2}\right\}\geq \dfrac{\hat\eps_*^2}{4C_3}=:4\bar q, \, \text{ for } (s,x)\in\h_2^\gamma,
$$
which implies that
\begin{equation}\label{lm3.5-e11}
\PP^{v_\gamma}_{s,x}\{\eta\leq T^*\}\geq 4\bar q, \text{ for }\, (s,x)\in\h^\gamma_2,
\end{equation}
where
$$\eta=\inf\left\{t\geq 0: X^{2}(t)\geq \frac{\hat\eps_*}{2}\right\}.$$
On the other hand, we have from It\^o's formula that
\begin{equation}\label{lm3.5-e12}
\begin{aligned}
\E^{v_\gamma}_{s,x}S(\eta\wedge T^*)=&s+\E^{v_\gamma}_{s,x}\int_0^{\eta\wedge T^*} \left[\dfrac{S_0-S(t)}\theta-X(t)f_1(S(t),X(t))\right]dt\\
\leq& s+\dfrac{S_0}\theta\E^{v_\gamma}_{s,x}(\eta\wedge T^*)\\
\leq& K_4+\dfrac{S_0T^*}{\theta}:=H_1.
\end{aligned}
\end{equation}
Defining
$$\tau^\eta_\h=\inf\{t\geq\eta\wedge T^*: (S(t), X(t))\in\h\},$$
we have
\begin{equation}\label{lm3.5-e13}
\begin{aligned}
\E^{v_\gamma}_{s,x}\tau^\eta_\h=&\E^{v_\gamma}_{s,x}(\eta\wedge T^*)+\E^{v_\gamma}_{s,x}(\tau^\eta_\h-(\eta\wedge T^*))\\
\leq &T^*+\E^{v_\gamma}_{s,x}\left[\E_{S(\eta\wedge T^*),X(\eta\wedge T^*)}^{v_\gamma}\tau_\h\right]\\
\leq &T^*+\E^{v_\gamma}_{s,x}\left[2\kappa_0S(\eta\wedge T^*)+X(\eta\wedge T^*)+C_2\right]\\
\leq &T^*+2\kappa_0H_1+1+C_2.
\end{aligned}
\end{equation}
It follows from \eqref{lm3.5-e11}, \eqref{lm3.5-e12}, and Markov's inequality that
\begin{equation}\label{lm3.5-e14}
\PP^{v_\gamma}_{s,x}\left\{\eta\leq T^*\,\text{ and } S(\eta)\leq \frac{H_1}{2\bar q}\,\right\}\geq 2\bar q, \text{ for }\, (s,x)\in\h^\gamma_2.
\end{equation}
By Lemma \ref{lm3.6}, there exists a $\delta_*>0$ depending only on $\frac{H_1}{2\bar q}$ and $\hat\eps_*$ such that
\begin{equation}\label{lm3.5-e15}
\PP_{s,x}^{v_\gamma}\{X(\tau_\h)\geq\delta_*\}>\frac12, \text{ for } s\leq\frac{H_1}{2\bar q}, x^2\geq \frac{\hat\eps_*}{2}, \gamma\in(0,1).
\end{equation}
We can choose $\delta_*<x_*$ and define $\h_3=\{(s,x)\in\h: x\geq \delta_*\}$.
Then $\h_3$ is a compact subset of $\R^{2,\circ}_+$.
We have from the strong Markov property of $(S(t), X(t))$ under Markov control $v_\gamma$
and from \eqref{lm3.5-e14} and \eqref{lm3.5-e15}
that
\begin{equation}\label{lm3.5-e16}
\begin{aligned}
\PP^{v_\gamma}_{s,x}\left\{(S(\tau^\eta_\h), X(\tau^\eta_\h))\in\h_3\right\}
\geq& \E^{v_\gamma}_{s,x} \left[\1_{\left\{\eta\leq T^*\,\text{ and } S(\eta)\leq \frac{H_1}{2\bar q}\,\right\}}\1_{\{(S(\tau^\eta_\h), X(\tau^\eta_\h))\in\h_3\}}\right]\\
=&\E^{v_\gamma}_{s,x}\left[\1_{\left\{\eta\leq T^*\,\text{ and } S(\eta)\leq \frac{H_1}{2\bar q}\,\right\}}\E^{v_\gamma}_{S(\eta),X(\eta)}\1_{\{X(\tau_\h)\geq\delta_*\}}\right]\\
\geq&\dfrac12\E^{v_\gamma}_{s,x}\left[\1_{\left\{\eta\leq T^*\,\text{ and } S(\eta)\leq \frac{H_1}{2\bar q}\,\right\}}\right]\\
\geq&\bar q.
\end{aligned}
\end{equation}
With $A=\left\{(S(\tau^\eta_\h), X(\tau^\eta_\h))\in\h_3\right\}$
and $A^c=\Omega\setminus A$, we have the estimate
$$
\begin{aligned}
\bar V_\gamma(s,x)=&\E^{v_\gamma}_{s,x}\int_0^{\tau^{\eta}_\h} e^{-\gamma t}S(t)dt + \E^{v_\gamma}_{s,x} \left[e^{-\gamma \tau^{\eta}_\h} V_\gamma\left(S(\tau^{\eta}_\h), X(\tau^{\eta}_\h)\right)\right]-V_\gamma(s_*,x_*)\\
\geq&\E^{v_\gamma}_{s,x} \left[\1_Ae^{-\gamma \tau^{\eta}_\h} \inf_{(s',x')\in\h_3} V_\gamma(s',x')\right]
+ \E^{v_\gamma}_{s,x} \left[\1_{A^c}e^{-\gamma \tau^{\eta}_\h} \inf_{(s',x')\in\h} V_\gamma(s',x')\right]
-V_\gamma(s_*,x_*)\\
=&\E^{v_\gamma}_{s,x} \left[\1_A\left(\inf_{(s',x')\in\h_3} V_\gamma(s',x')-V_\gamma(s_*,x_*)\right)\right]\\
& +\E^{v_\gamma}_{s,x} \left[\1_{A^c}\left(\inf_{(s',x')\in\h} V_\gamma(s',x')-V_\gamma(s_*,x_*)\right)\right]\\
&-\E^{v_\gamma}_{s,x} \left[\1_{A}(1-e^{-\gamma \tau^{\eta}_\h}) \inf_{(s',x')\in\h_3} V_\gamma(s',x')\right]\\
&- \E^{v_\gamma}_{s,x} \left[\1_{A^c}(1-e^{-\gamma \tau^{\eta}_{\h}}) \inf_{(s',x')\in\h} V_\gamma(s',x')\right].
\end{aligned}
$$
Since $\bar V_\gamma(s,x)= V_\gamma(s,x)-V_\gamma(s_*, x_*)\to V(s,x)$ as $\gamma\to0$
uniformly in each compact set,
there exists  an $H_2>0$ such that $|\bar V_\gamma(s,x)|<H_2$ for $(s,x)\in\h_3$ when $\gamma$ is sufficiently small.
We also have
$0\leq\inf_{(s',x')\in\h_3} V_\gamma(s',x')\leq V(s_*,x_*)\leq\frac{\bar\rho}{\gamma}$
when $\gamma$ is sufficiently small.
This together with \eqref{lm3.5-e13}, \eqref{lm3.5-e16}, and
$\inf_{(s,x)\in\h} \bar V_\gamma(s,x)\leq 0$
yields that
$$
\begin{aligned}
\bar V_\gamma(s,x)\geq& -H_2\PP^{v_\gamma}_{s,x}(A)+\inf_{(s',x')\in\h} \bar V_\gamma(s',x') \PP^{v_\gamma}_{s,x}(A^c)
- \bar\rho\E^{v_\gamma}_{s',x'} \left[\dfrac{1-e^{-\gamma \tau^{\eta}_\h}}\gamma\right]\\
\geq& -H_2\PP^{v_\gamma}_{s,x}(A)+\inf_{(s',x')\in\h} \bar V_\gamma(s',x') \PP^{v_\gamma}_{s,x}(A^c)
- \bar\rho\E^{v_\gamma}_{s,x} \tau^{\eta}_\h\\
\geq& -H_2-\bar\rho\left (T^*+2\kappa_0H_1+1+C_2\right)+\inf_{(s',x')\in\h} \bar V_\gamma(s',x') \PP^{v_\gamma}_{s,x}(A^c), \\
\geq& -H_2-\bar\rho\left (T^*+2\kappa_0H_1+1+C_2\right)+\inf_{(s',x')\in\h} \bar V_\gamma(s',x') (1-\bar q), \text{ for any } (s,x)\in\h_2^\gamma,\\
\end{aligned}
$$
which combined with \eqref{lm3.5-e9} leads to
$$\inf_{(s,x)\in\h} \bar V_\gamma(s,x) \geq -H_2-\bar\rho\left (T^*+2\kappa_0H_1+1+C_2\right)+\inf_{(s,x)\in\h} \bar V_\gamma(s,x) (1-\bar q),$$
or
\begin{equation}\label{lm3.5-e17}
\begin{aligned}
\inf_{(s,x)\in\h} \bar V_\gamma(s,x)\geq&  - \dfrac1{\bar q}\left(H_2+\bar\rho\left (T^*+2\kappa_0H_1+1+C_2\right)\right)=:-H_3.
\end{aligned}
\end{equation}

Let $v_c\equiv M$ be the constant control.
Similar to the proof of Lemma \ref{lm3.3}, we can show that
$$H_4:=\sup_{(s,x)\in\h}\E^{v_c}_{s,x}\tau_{\h_3}<\infty.$$
We have
$$\E^{v_c}_{s,x} S(\tau_{\h_3})\leq s+\dfrac{S_0}\theta\E^{v_c}_{s,x}\tau_{\h_3}-\frac1\theta\E^{v_c}_{s,x} \int_0^{\tau_{\h_3}} S(t),
$$
which implies
\begin{equation}\label{lm3.5-e19}
\E^{v_c}_{s,x}\int_0^{\tau_{\h_3}} S(t)\leq \theta s+ S_0\E^{v_c}_{s,x}\tau_{\h_3}
\leq \theta K_4+S_0 H_4, (s,x)\in\h.
\end{equation}
By \cite[Eq. (3.7.47)]{ABG},
\begin{equation}\label{lm3.5-e18}
\begin{aligned}
V(s,x)\leq& \E^{v_c}_{s,x}\left[\int_0^{\tau_{\h_3}}S(t)dt+V(S(\tau_{\h_3}), X(\tau_{\h_3}))\right]\\
\leq& \theta K_4+S_0 H_4+\sup_{(s',x')\in\h_3} V(s',x').
\end{aligned}
\end{equation}
From \eqref{lm3.5-e17} and \eqref{lm3.5-e18},
we obtain that
$$\sup_{(s,x)\in\h}|V(s,x)|<\infty.$$
\end{proof}

\begin{proof}[Proof of Theorem \ref{HJB}]
For any $(s,x)\in\R^{2,\circ}_+$,
we have
\begin{equation}\label{HJB-e4}
\begin{aligned}
\bar V_\gamma(s,x)=&\E^{v_\gamma}_{s,x}\int_0^{\tau_\h} e^{-\gamma t}S(t)dt + \E^{v_\gamma}_{s,x} \left[e^{-\gamma \tau_\h} V_\gamma\left(S(\tau_\h), X(\tau^T_\h)\right)\right]-V_\gamma(s_*,x_*)\\
\geq&\inf_{(s',x')\in\h} \bar V_\gamma(s',x')
-\E^{v_\gamma}_{s,x}\frac{\rho}{\gamma} (1-e^{-\gamma\tau_\h})\\
\geq&\inf_{(s',x')\in\h} \bar V_\gamma(s',x')
-\E^{v_\gamma}_{s,x}\tau_\h\\
\geq&-(H_3 + 2\kappa_0 s+x+C_1)\quad\text{due to \eqref{lm3.5-e17} and \eqref{lm3.3-e0}.}
\end{aligned}
\end{equation}
As a result,
$$V(s,x)\geq -(H_3 +2\kappa_0 s+x+C_1+1).$$
Similar to \eqref{lm3.5-e19} and \eqref{lm3.5-e18},
we have
\begin{equation}\label{HJB-e5}
\begin{aligned}
V(s,x)\leq &\E^{v_c}_{s,x}\left[\int_0^{\tau_{\h}}S(t)dt+V(S(\tau_{\h}), X(\tau_{\h}))\right]\\
\leq&\sup_{(s',x')\in\h}|V(s',x')|+ \theta s+ S_0\E^{v_c}_{s,x}\tau_{\h}\\
\leq & \sup_{(s',x')\in\h}|V(s',x')|+ \theta s+ S_0(2\kappa_0s+x+\hat C_1).
\end{aligned}
\end{equation}
In view of \eqref{HJB-e4} and \eqref{HJB-e5},
we have $V(s,x)\in\mathcal{C}_{1+p}$.
Then,
we can use
arguments similar to \cite[Theorem 3.7.11 and Theorem 3.7.12]{ABG}
and Lemma \ref{lm3.2} to obtain the desired result.
\end{proof}

\section{Concluding Remarks}\label{sec:con}
To validate and to improve  model \eqref{ww2},
verification using real data is needed.
To verify the model, the parameters of the system need to be estimated first.
A statistical estimator can be constructed.
To estimate the parameters using real data, we observe the solutions of \eqref{ww2} in discrete epoch, and carry out the estimation accordingly.
 That is, view the observation (the real data) as solution of \eqref{ww2}, then use
the explicit Euler method to discretize the diffusion process \eqref{ww2}, and utilize
for example,
the
maximum likelihood method to estimate the parameter.
An alternative approach is to use the generalized method of moments.

The simplified model \eqref{ww2} may not be sufficient to perceive the complicated process of wastewater treatment.
Considering more complex models
renders better understanding but also poses more challenges.

In Section \ref{sec:control}, we worked with controlled diffusions without switching for notational simplicity.
The proofs carry over if one
considers the controlled switching diffusion counterpart of \eqref{ww-c1}.
We have proved the existence and uniqueness of solutions to the HJB equation for \eqref{ww-c1}. For the controlled switching diffusions, we need to deal with a system of HJB equations.
While the optimal Markov control can be obtained theoretically using the system of HJB equations, it is quite difficult to
find a closed-form solution explicitly.
Constructing a numerical scheme is a viable alternative.

\appendix
\section{Proofs of Theorem \ref{thm2.1}}
\begin{proof}[Proof of Theorem \ref{thm2.1}]
Since the coefficients of \eqref{ww3} are locally Lipschitz continuous in $(s,x,i)\in \R_+^2\times\M$,
the system given by \eqref{ww3} and \eqref{eq:tran} has a unique continuous solution up to explosion time $\tau_e$, where $\tau_e=\inf\{t\geq 0: S(t)\vee X(t)=\infty\}$, with the convention $\inf\emptyset=\infty$. The solution is also a strong Markov process; see \cite{XM, ZY}. If we define
$$\tau_k=\inf\Big\{t\geq 0: S(t)\vee X(t)>k\Big\},$$
then $\tau_e=\lim_{k\to\infty}\tau_k$.
Consider $\hat V_1(s,x,i)=\kappa_0s+x$ then, by the generalized It\^o formula, we have
$$
\begin{aligned}
\Lom \hat V_1(s,x,i)&=\dfrac{S_0\kappa_0}{\theta}-\dfrac{\kappa_0}{\theta}f_0(s,i)-\wdt k_d(i)x+x(f_2(s,x,i)-\kappa_0f_1(s,x,i))
\\&\leq \dfrac{S_0\kappa_0}{\theta}.
\end{aligned}
$$
Hence,
$$
\begin{aligned}
\E_{s,x,i} \hat V_1\big(S(\tau_k\wedge t), X(\tau_k\wedge t),\alpha(\tau_k\wedge t)\big)
\leq \hat V_1(s,x,i)+\dfrac{S_0\kappa_0t}{\theta},
\end{aligned}
$$
which implies that
$$
\begin{aligned}
\PP_{s,x,i}\big\{\tau_k<t\big\}&\leq \PP_{s,x,i}\Big\{\hat V_1\big(S(\tau_k\wedge t), X(\tau_k\wedge t),\alpha(\tau_k\wedge t)\big)\geq k\Big\}
\\&\leq \dfrac{\hat V_1(s,x,i)+\frac{S_0\kappa_0t}{\theta}}{k}\to 0\;\text{as}\;k\to\infty.
\end{aligned}
$$
Therefore, we have $\PP_{s,x,i}\{\tau_e\leq t\}=0$ or $\PP_{s,x,i}\{\tau_e>t\}=1\;\forall t>0$. As a consequence, $\PP_{s,x,i}\{\tau_e=\infty\}=1$. Hence, the system given by \eqref{ww3} and \eqref{eq:tran} has a unique global continuous solution.

Now, we move to the part of positivity of solutions. First, suppose that $s,x>0$.
For any $n\in\N$, we define the following truncated functions
$$f_0^{(n)}(s,i)=f_0(s\wedge n,i)\;;\;f_{1}^{(n)}(s,x,i)=f_1(s\wedge n,x\wedge n,i)\;;\;f_{2}^{(n)}(s,x,i)=f_2(s\wedge n,x\wedge n,i),$$
and let $\big(S^{(n)}(t),X^{(n)}(t)\big)$ be the solution of equation \eqref{ww3} and \eqref{eq:tran} with $f_0,f_1,f_2$ replaced by $f_0^{(n)},f_1^{(n)},f_2^{(n)}$, respectively.
Denote
$$\eta^{(n)}=\inf\big\{t\geq 0: S^{(n)}(t)\wedge X^{(n)}(t)\leq 0\big\},$$
$$\eta_k^{(n)}=\inf\Big\{t\geq 0: S^{(n)}(t)\wedge X^{(n)}(t)<\dfrac 1k\Big\}.$$
Then $\eta^{(n)}=\lim_{k\to\infty}\eta_k^{(n)}$. Consider
$$\hat V_2^{(n)}(s,x,i)=s-c_1^{(n)}-c_1^{(n)}\ln \dfrac s{c_1^{(n)}}+c_2(x-1-\ln x),$$
where
$$c_2=\dfrac 1{\kappa_0}\;\text{and}\;c_1^{(n)}=\frac{c_2\check k_d}{\max_{i\in\M}\sup_{s>0}\frac{f_1^{(n)}(s,x,i)}{s}}.$$
Let
$$c_2^{(n)}:=\max_{i\in\M}\sup_{s>0}\frac{f_0^{(n)}(s,i)}{s}.$$
 Since $f_1^{(n)}(s,x,i)$ is global Lipschitz and $f_1^{(n)}(0,x,i)=0$, it is readily seen that
 $\displaystyle\sup_{s>0}\frac{f_1^{(n)}(s,x,i)}{s}<\infty$.
Similarly, we also obtain  $\sup_{s>0}\frac{f_0^{(n)}(s,i)}{s}<\infty$, so
the above constants are well-defined.

 By the generalized It\^o formula, we have
$$
\begin{aligned}
\Lom^{(n)} \hat V_2^{(n)}(s,x,i)=&\Big(1-\dfrac {c_1^{(n)}}s\Big)\left(\dfrac {S_0}\theta-\dfrac{f_0^{(n)}(s,i)}\theta-xf_1^{(n)}(s, x,i)\right)+\dfrac{c_1^{(n)}\sigma_1^2(i)}{2}\\
&+\Big(c_2-\dfrac{c_2}{x}\Big)\left(xf_2^{(n)}(s, x,i)-x\wdt k_d(i) \right)+\dfrac{c_2\sigma_2^2(i)}{2}
\\
\leq& \dfrac{S_0+c_1^{(n)}c_2^{(n)}}{\theta}+c_2\hat k_d+\dfrac{c_1^{(n)}\hat\sigma_1^2+c_2\hat \sigma_2^2}{2}\\
&+x\Big(c_2f_2^{(n)}(s,x,i)-f_1^{(n)}(s,x,i)\Big)+x\Big(\dfrac{c_1^{(n)}f_1^{(n)}(s,x,i)}{s}-c_2k_d(i)\Big) \\
\leq& \dfrac{S_0+c_1^{(n)}c_2^{(n)}}{\theta}+c_2\hat k_d+\dfrac{c_1^{(n)}\hat\sigma_1^2+c_2\hat \sigma_2^2}{2}:=K^{(n)} \text{ for } s,x>0,
\end{aligned}
,$$
where the operator $\Lom^{(n)}$ is defined as $\Lom$ with $f_0$, $f_1$, and $f_2$ replaced by $f_0^{(n)}$, $f_1^{(n)}$, and $f_2^{(n)}$, respectively.
Applying It\^o's formula again, we have
$$
\begin{aligned}
\E_{s,x,i} &\hat V_2^{(n)}\big(S^{(n)}(\eta_k^{(n)}\wedge t), X^{(n)}(\eta_k^{(n)}\wedge t),\alpha(\eta_k^{(n)}\wedge t)\big)\\
=&\hat V_2^{(n)}(s,x,i)+\E_{s,x,i}\int_0^{\eta_k^{(n)}\wedge t}\Lom \hat V_2^{(n)}\big(S^{(n)}(u), X^{(n)}(u),\alpha(u)\big)du\\
\leq& \hat V_2^{(n)}(s,x,i)+K^{(n)}t.
\end{aligned}
$$
Since the definition of $\hat V_2^{(n)}$, we observe that if $\eta_k^{(n)}<t$, then
$$
\hat V_2^{(n)}\big(S^{(n)}(\eta_k^{(n)}\wedge t), X^{(n)}(\eta_k^{(n)}\wedge t),\alpha(\eta_k^{(n)}\wedge t)\big)\geq \big(c_1^{(n)}\ln kc_1^{(n)}- c_1^{(n)}\big)\wedge \big(c_2\ln k-c_2\big) .
$$
Hence, we obtain that
\begin{equation*}
\begin{aligned}
\PP_{s,x,i}\{\eta_k^{(n)}<t\}\leq \dfrac{\hat V_2^{(n)}(s,x,i)+K^{(n)}t}{\big(c_1^{(n)}\ln kc_1^{(n)}- c_1^{(n)}\big)\wedge \big(c_2\ln k-c_2\big) }\to 0\;\text{as}\;k\to \infty.
\end{aligned}
\end{equation*}
As a result, for any $n\in \N$, $\PP_{s,x,i}\{\eta^{(n)}_\infty=\infty\}=1$.
That is,
$$\PP_{s,x,i}\big\{S^{(n)}(t),X^{(n)}(t)>0:\forall t>0\big\}=1\text{ for any }n\in\N,$$
and hence
$$\PP_{s,x,i}\big\{S^{(n)}(t),X^{(n)}(t)>0:\forall t>0,n\in\N\big\}=1.$$
Now, for any $t>0$ and
$$\omega\in \{\tau_e=\infty\}\cap\big\{S^{(n)}(t),X^{(n)}(t)>0:\forall t>0,n\in\N\big\},$$
there exists $n_0=n_0(\omega,t)$ such that
$$S(y)(\omega)\vee X(y)(\omega)<n_0\;\forall 0\leq y\leq t.$$
As a consequence, $S(t)(\omega)=S^{(n_0)}(t)(\omega)>0$ and $X(t)(\omega)=X^{(n_0)}(t)(\omega)>0$. This, combined with
$\PP_{s,x,i}\{\tau_e=\infty\}=1$, implies that
\begin{equation}\label{e2-thm2.1}
\PP_{s,x,i}\big\{S(t)>0: t>0 \big\}=\PP_{s,x,i}\big\{X(t)>0: t>0 \big\}=1\;\forall s,x>0.
\end{equation}
If $s>0,x=0$, the result $\PP_{s,x,i}\big\{S(t)>0: t>0 \big\}=1$ is similar proved by choosing $c_2=0$. Moreover, it is obvious that $\PP_{s,x,i}\big\{X(t)=0: t>0 \big\}=1$.

Consider the case when the initial value $s=0$ and $x\geq 0$.
Let $\eps>0$ be sufficiently small such that
\begin{equation}\label{e3-thm2.1}\dfrac{S_0-f_0(\tilde s,i)}\theta-\tilde xf_1(\tilde s,\tilde x,i)\geq \frac{S_0}{2\theta},
\end{equation}
for any $(\tilde s,\tilde x,\tilde i)\in\R^2\times\M$ satisfying
$\tilde s+|\tilde x-x|<\eps$.
Let $$\tilde\tau_1=\inf\{t>0:S(t)+|X(t)-x|\geq\eps\}.$$
By the continuity of $(S(t), X(t))$, $\PP_{0,x,i}\{\tilde\tau_1>0\}=1$.
Using the variation of constants formula (see \cite[Chapter 3]{XM}),
we can write $S(t)$ in the form
\begin{equation}\label{e4-thm2.1}\barray
\ad
S(t)=\Phi(t)\left[\int_0^t\Phi^{-1}(u)\left(\dfrac{S_0-f_0\big(S(u),
\alpha(u)\big)}\theta-X(u)f_1\big(S(u),X(u),\alpha(u)\big)\right)du\right] \\
\aad \ \hfill \text{ for } t\in[0,\tilde\tau_1),\earray
\end{equation}
where
$\Phi(t)=\exp\left(-\int_0^t\dfrac{\sigma_1^2(\alpha(u))}2du+\int_0^t\sigma_1(\alpha(u))dW_2(u)\right).$
It follows
from \eqref{e3-thm2.1} that
$$\frac{S_0-f_0\big(S(u),\alpha(u)\big)}\theta-X(u)f_1\big(S(u),X(u),\alpha(u)\big)>0\;\text{if}\;t\in(0,\tilde\tau_1].$$
 This and \eqref{e4-thm2.1} imply that
$$\PP_{0,x,i}\{S(t)>0, t\in(0,\tilde\tau_1]\}=1,$$
which combined with \eqref{e2-thm2.1}
and the strong Markov property of $(S(t),X(t),\alpha(t))$
yields that
\begin{equation*}
\PP_{0,x,i}\{S(t)>0, t\in(0,\infty)\}=1.
\end{equation*}
The theorem is therefore proved.
\end{proof}

\section{Proof of Lemmas in Section \ref{sec:control}}
\begin{proof}[Proof of Lemma \ref{lm3.3}]
Since $f_1(0,x)=0$, we can choose $\delta_1>0$ such that
 $s+\theta xf_1(s,x)\leq \frac{S_0}2$ for $s\leq\delta_1, x\leq K_4$, where $K_4$ is defined in \eqref{defh0}.
Due to the uniform bound from \eqref{lm3.1-e0} and
the stochastic continuity of $S(t)$, (which can be seen to be uniform
in $v$
since the equation of $S(t)$ does not depend directly on $v$),
there exist $t_1\in(0,\frac12)$ and $\hat \delta_1\in(0,\delta_1)$ such that
\begin{equation}\label{lm3.3-e1a}
\PP^v_{s,x}\left\{S(t)<\delta_1, X(t)\leq K_4+1, t\in[0,t_1]\right\}\geq\frac34, (s,x)\in\h_0,s<\hat\delta_1.
\end{equation}
 Let $K_5$ be sufficiently large such that
\begin{equation}\label{lm3.3-e2a}
\PP\left\{\Phi(t)\vee \Phi^{-1}(t)<K_5, t\in [0,t_1]\right\}\geq \frac{3}{4},
\end{equation}
where in this section,
$\Phi(t)=\exp\left(-\frac{\sigma_1^2}2t+\int_0^t\sigma_1dW_1(u)\right).$
If $\Phi(t)\vee \Phi^{-1}(t)<K_5$ and $S(t)+\theta X(t)f_1(S(t),X(t))\leq \frac{S_0}2$ for $t\in[0,t_1]$,
we have
\begin{equation}\label{lm3.3-e3a}
S(t_1)=\Phi(t)\int_0^t\Phi^{-1}(t)\left(\frac{S(u)-S_0}\theta-X(u)f_1(S(u),X(u))\right)du\geq \frac{t_1S_0}{2\theta K_5^2}.
\end{equation}
Let $\delta_0:=\frac{t_1S_0}{2\theta K_5^2}\wedge \hat\delta_1$ and
$\h=\{(s,x): 2\kappa_0s+x\leq K_4, s\geq\delta_0\}$.
Since
the diffusion is nondegenerate on $\R^{2,\circ}_+$, it is well-known (see e.g., \cite[Lemma 2.6.5]{ABG}) that
there exist $t_2>0$ and $p_1>0$ such that	
$$\PP^v_{s,x}\left\{\tau_\h<t_2\right\}>2p_1, \text{ provided } \delta_0 \leq s\leq\delta_1, K_4-2\kappa_0\delta_1\leq x\leq K_4+1.$$
Because of the definition of $\h$,
$$\PP^v_{s,x}\left\{\tau_\h<t_2\right\}>2p_1, \text{ provided } \delta_0 \leq s\leq\delta_1, x\leq K_4+1.$$
This together with \eqref{lm3.3-e1a}, \eqref{lm3.3-e2a}, \eqref{lm3.3-e3a},
 and the definition of $\delta_0$  implies
\begin{equation}\label{lm3.3-e8'}
\PP_{s,x}^v\left\{\tau_\h<t_1+t_2\right\}\geq p_1, \, (s,x)\in \h_0.
\end{equation}

Define stopping times
$$
\begin{aligned}
&\tau^{(1)}_{\h_0}=\inf\{t\geq t_1+t_2: (S(t), X(t))\in \h_0\},\\
&\tau^{(n+1)}_{\h_0}=\inf\{t\geq \tau^{n}_H+t_1+t_2: (S(t), X(t))\in \h_0\}.
\end{aligned}
$$
In view of
\eqref{lm3.1-e0}
$$\E_{s,x}^v \big(2\kappa_0S(t_1+t_2)+X(t_1+t_2)\big)\leq (2\kappa_0 s+x)+K_1(t_1+t_2)\leq K_4+K_1(t_1+t_2),$$
which together with \eqref{e3.2} implies
$$\E_{s,x}^v \tau^{(1)}_{\h_0}\leq K_4+K_1(t_1+t_2), (s,x)\in\h_0.$$
Then
the strong Markov property implies that
\begin{equation}\label{lm3.3-e9}
\E_{s,x}^v \left[\tau^{(n+1)}_{\h_0}-\tau^{(n)}_{\h_0}-1\Big|\F_{\tau^{(n)}_{\h_0}+1}\right]\leq K_4+K_1(t_1+t_2), (s,x)\in\h_0, n\geq1.
\end{equation}
Let
$$B_k=\left\{\big(S(t), I(t)\big)\in\h  \text{ for some } t\in\left[\tau^{(k)}_{\h_0},\tau^{(k)}_{\h_0}+1\right)\right\}, k\geq0.$$
Let $\sigma^{(1)}=\tau^{(1)}_{\h_0}, \sigma^{(k)}=\tau^{(k)}_{\h_0}-\tau^{(k-1)}_{\h_0}-1, k\geq 2$.
Note that
if $B_k$ occurs, then $\tau_\h\leq \tau^{(k)}_{\h_0}+1$.
By the strong Markov property of $\big(S(t), X(t)\big)$
$$
\begin{aligned}
\disp \PP_{s,x}^v\{\cap_{k=1}^{n+1}B^c_k\}
=&\disp \E_{s,x}^v \left[\1_{\{\cap_{k=1}^nB^c_k\}}\E\left[B^c_{n+1}\big|\F_{\tau^{(n+1)}}\right]\right]\\
\leq&\disp p_1\E_{s,x}^v \left[\1_{\{\cap_{k=1}^nB^c_k\}}\right]\\
=&\disp p_1\PP_{s,x}^v\{\cap_{k=1}^{n}B^c_k\}.
\end{aligned}
$$
By induction, we have
$$
\PP_{s,x}^v\{\cap_{k=1}^{n+1}B^c_k\}\leq p_1^{n-1},
$$
which leads to
$$
\PP_{s,x}^v\{\cap_{k=1}^{\infty}B^c_k\}=0.
$$
As a result,
\begin{equation}\label{e13-thm3.1}
\sum_{n=1}^\infty\PP_{s,x}^v\{B_n\cap_{k=1}^{n-1}B^c_k\}
=1,
\end{equation}
and
\begin{equation}\label{e14-thm3.1}
\begin{aligned}
\E_{s,x}^v \tau_\h
=&\sum_{n=1}^\infty\E_{s,x}^v \left[\tau_\h\1_{\{B_n\cap_{k=1}^{n-1}
B^c_k\}}\right]\\
\leq& \sum_{n=1}^\infty\E_{s,x}^v \left[[\tau^{(n)}+1]\1_{\{B_n\cap_{k=1}^{n-1}B^c_k\}}\right]\\
\leq& 1+\sum_{n=1}^\infty\E_{s,x}^v \left[\sum_{l=1}^n\sigma^{(l)}\1_{\{B_n\cap_{k=1}^{n-1}B^c_k\}}\right]\\
=&1+\sum_{l=1}^\infty\E_{s,x}^v \left[\sigma^{(l)}\sum_{n=l}^\infty\1_{\{B_n\cap_{k=1}^{n-1}B^c_k\}}\right]\\
=&1+\sum_{l=1}^\infty\E_{s,x}^v \left[\sigma^{(l)}
\1_{\{\cap_{k=1}^{l-1}B^c_k\}}\right]\,\text{ (due to \eqref{e13-thm3.1})}.
\end{aligned}
\end{equation}
In view of \eqref{lm3.3-e9} we have
\begin{equation}\label{e16-thm3.1}
\begin{aligned}
\E_{s,x}^v \left[\sigma^{(l)}\1_{\{\cap_{k=1}^{l-1}B^c_k\}}\right]
=&\E_{s,x}^v \left[\1_{\{\cap_{k=1}^{l-1}B^c_k\}}\E\left[\sigma^{(l)}\big|\F_{\tau^{l-1}_{\h_0}}\right]\right]\\
\leq& \left( K_4+K_1(t_1+t_2)\right)\E_{s,x}^v \left[\1_{\{\cap_{k=1}^{l-1}B^c_k\}}\right]\\
\leq&\left(K_4+K_1(t_1+t_2)\right){p_1}^{l-1},\,l\geq 1.
\end{aligned}
\end{equation}
Therefore, it follows from
\eqref{e14-thm3.1} and \eqref{e16-thm3.1} that
$$
\begin{aligned}
\E_{s,x}^v \tau_\h &\!\disp \leq 1+\sum_{l=1}^\infty\E_{s,x}^v \left[\sigma^{(l)}
\1_{\{\cap_{k=1}^{l-1}B^c_k\}}\right]\\
&\! \disp
\leq 1+\sum_{l=1}^\infty \left( K_4+K_1(t_1+t_2)\right){p_1}^{l-1}:=C_1<\infty.
\end{aligned}
$$	
This and \eqref{e3.2} imply \eqref{lm3.3-e0}.
Under the constant control $v^c\equiv M$,
similar to that of $S(t)$, the drift of $X(t)$ is positive when $X(t)$ is small.
As a result, a similar argument can be deployed to obtain \eqref{lm3.3-e0a}.
\end{proof}

\begin{proof}[Proof of Lemma \ref{lm3.6}]
We have
$$
\begin{aligned}
\op^u\left(\ln\frac{1+x}x\right)=&
-\dfrac1{x(1+x)}\left[u+xf_2(s,x)-x\wdt k_d\right]-\dfrac{\sigma_2^2x^2}{(1+x)^2}
+\dfrac{\sigma_2^2}2\\
\leq&\wdt k_d+\dfrac{\sigma_2^2}2.
\end{aligned}
$$
Then by Dynkin's formula and some standard calculations,
we can show that
$$
\begin{aligned}
\E^v_{s,\eps}\ln \frac{1+X(\tau_\h)}{X(\tau_\h)}
\leq& \ln\frac{1+x}{x}+ \left(\wdt k_d+\dfrac{\sigma_2^2}2\right)\E^v_{s,\eps}\tau_\h\\
\leq& \ln \frac{1+\eps}\eps +\left(\wdt k_d+\dfrac{\sigma_2^2}2\right)\sup_{s'\leq L_1}\E^v_{s',\eps}\tau_\h\\
=:&L_2<\infty,
\end{aligned}
$$
for $s\leq L_1$.
By Markov inequality, we have
$$\PP^v_{s,\eps}\left\{\ln \frac{1+X(\tau_\h)}{X(\tau_\h)}\leq 2L_2\right\}\geq\frac12,$$
which implies \eqref{lm3.6-e0}.
\end{proof}

\begin{proof}[Proof of Lemma \ref{lm3.2}]
By Lemma \ref{lm3.1},
the family $\{\zeta_T^m, T>0\}$
is tight on $\PUB$
for any admissible relaxed control $m$.
As a result,
we can decompose any limit point $\hat\zeta\in\PUB$ as 	
$$\hat\zeta=\delta\zeta'+(1-\delta)\zeta'',$$
where $\zeta'\in\PU$
and $\zeta''\big((0,\infty)\times\{0\}\times[0,M]\big)=1$.
Following the arguments in \cite[Lemma 3.4.6]{ABG},
we can show that $\zeta'\in\mathcal{G}$.
Because of \eqref{Sinfty} and the uniform boundedness \eqref{lm3.1-e1} and \eqref{lm3.1-e3}, we have
$s+xf_1(s,x)$ is $\hat\zeta$-integrable and
\begin{equation}\label{LSzetahat}
\int_{\R^{2}_+\times[0,M]}\left(\dfrac{S_0-s'}\theta-xf_1(s,x)\right)\hat\zeta(ds',dx',du)=0.
\end{equation}
Since $\zeta'$ is in $\mathcal{G}$, we have from \eqref{LSeq} that
\begin{equation}\label{LSeqzeta}
\int_{\R^{2,\circ}_+\times[0,M]}\left(\dfrac{S_0-s'}\theta-xf_1(s,x)\right)\zeta'(ds',dx',du)=0.
\end{equation}
As a result,
$$\int_{(0,\infty)\times\{0\}\times[0,M]}\dfrac{S_0-s'}
\theta\zeta''(ds',dx',du)=\int_{\R^{2}_+\times[0,M]}
\left(\dfrac{S_0-s'}\theta-xf_1(s,x)\right)\zeta''(ds',dx',du)=0,
$$
or equivalently,
$$\int_{(0,\infty)\times\{0\}\times[0,M]}s'\zeta''(ds',dx',du)=S_0.$$

To prove \eqref{lm3.2-e2},
note that we can find a sequence $t_n\uparrow \infty$ satisfying
$$
\liminf_{T\to\infty} \dfrac1T\E_{s,x}^m\int_0^TS(t)dt=\lim_{n\to\infty} \dfrac1{t_n}\E_{s,x}^m\int_0^{t_n}S(t)dt $$
and
$\zeta^m_{t_n}$ converges weakly to a probability measure $\hat\zeta$, which can be decomposed as
\eqref{lm3.2-e1}.
By the weak convergence,  the uniform boundedness of $\E^m_{s,x} (2\kappa_0S(t)+X(t))^{1+p}$ in \eqref{lm3.1-e2}, and using \eqref{rhos0},
we obtain
$$
\begin{aligned}
\lim_{n\to\infty} \dfrac1{t_n}\E_{s,x}^m\int_0^{t_n}S(t)dt
=&\delta\int_{\R^2_+\times[0,M]} s'\zeta'(ds', dx', du)+(1-\delta)\int_{\R^2_+\times[0,M]} s'\zeta''(ds', dx', du)\\
\geq&\delta\rho^*+S_0(1-\delta)\geq\rho^*.
\end{aligned}
$$
The proof is complete.
\end{proof}

\para{\bf Acknowledgments.}
We are very grateful to the editors and reviewers for evaluating our manuscript and for the constructive comments and suggestions,
which have lead to much improvement in the exposition of the paper.

\end{document}